\documentclass[11pt]{article}
\usepackage[english]{babel}
\usepackage{tikz}
\usetikzlibrary{automata,positioning}
\usetikzlibrary{shapes}
\usetikzlibrary{calc,arrows.meta,positioning}
\let\footnote=\endnote

\usepackage{times}
\usepackage{indentfirst}
\usepackage[title]{appendix}
\usepackage{amsmath,amstext,amssymb,amsfonts, amsthm}
\usepackage{enumitem}
\usepackage{graphicx,color}
\usepackage[size=scriptsize]{caption}
\usepackage[size=footnotesize]{subcaption}
\newcommand{\ei}{\end{itemize}}     
\newcommand{\bc}{\begin{center}}  
\newcommand{\ec}{\end{center}}     

\newcommand{\ls}[1]
   {\dimen0=\fontdimen6\the\font \lineskip=#1\dimen0
   \advance\lineskip.5\fontdimen5\the\font \advance\lineskip-\dimen0
   \lineskiplimit=.9\lineskip \baselineskip=\lineskip
   \advance\baselineskip\dimen0 \normallineskip\lineskip
   \normallineskiplimit\lineskiplimit \normalbaselineskip\baselineskip
   \ignorespaces }

\usepackage{xcolor}
\numberwithin{equation}{section}

\newtheorem{lemma}{Lemma}[section]
\newtheorem{theorem}[lemma]{Theorem}
\newtheorem{corollary}[lemma]{Corollary}
\newtheorem{definition}[lemma]{Definition}

\newtheorem{remark}[lemma]{Remark}
\newtheorem{Assumption}[lemma]{Assumption}

\def\Y{\mathbb{Y}}

\def\X{\mathbb{X}}
\def\E{\mathbb{E}}
\def\A{\mathbb{A}}
\def\H{\mathbb{H}}
\def\R{\mathbb{R}}
\def\Z{\mathbb{Z}}
\def\W{\mathbb{W}}

\def\G{\mathcal{G}}

\def\B{\mathbb{B}}

\def\Gr{{\rm Gr}}
\def\v^*{{\rm \mathcal{V}}}

\begin{document}



\title{Dynamically Augmented CVaR for MDPs}

\author{Eugene~A.~Feinberg\footnote{Corresponding author, Department of Applied Mathematics and
Statistics, Stony Brook University, Stony Brook, NY 11794-3600, USA, eugene.feinberg@stonybrook.edu, https://orcid.org/0000-0002-8263-0772} \and  Rui~Ding\footnote{Department of Applied Mathematics and
Statistics,
 Stony Brook University,
Stony Brook, NY 11794-3600, USA,
rrd051488@gmail.com
}}

\maketitle

\begin{abstract}%
This paper studies optimization of  Conditional Value-at-Risk (CVaR)
for Markov Decision Processes (MDPs) with finite state and action sets. It introduces the Dynamically
augmented CVaR (DCVaR) risk objective function and provides an algorithm for its optimization. This paper investigates
a specially defined Robust MDP (RMDP), in which the state space is augmented with the tail
risk level. This RMDP, which we call the Dynamically augmented RMDP (DRMDP), was introduced
to the literature for  calculations of optimal CVaR values by value iteration more than ten years ago, but, as was understood
later, these value iterations compute lower bounds of minimal static CVaRs. DCVaR is defined
as a time consistent version of the static CVaR, and it is a lower bound of the static CVaR. It also can
be considered as a dynamic version of the nested CVaR. This paper provides an algorithm constructing a policy optimizing DCVaR  of total discounted costs. The correctness of this algorithm is
proved by studying a special mass transfer problem.  The results on RMDPs needed for this paper are provided in the appendix.
\end{abstract}%


\section{Introduction}\label{s1}

Conditional Value-at-Risk (CVaR), also sometimes called Average Value-at-Risk (AVaR), is one of the
most important and popular risk measures. For dynamic problems, it is possible  to consider the so-called
static CVaR, when the CVaR value is defined for each policy, and the goal is to find a policy optimizing
this value. However, there are problems with this approach dealing with computational complexity of finding
optimal policies and time inconsistency of static CVaRs. The issue of time consistency of risk measures is
discussed in several works, including
\cite{KF, RS06a, RS06b, Sh09, Sh12, Sh21}.
 An alternative approach is
to consider the nested CVaR, which is defined via a Robust Markov Decision Pprocess (RMDP) whose state space coincides
with the original state space of the problem. Another approach, which is also often used in practice, is to
consider an RMDP, whose states are pairs consisting of the original state augmented with the risk level. In
other words, the state is dynamically augmented by a number between 0 and 1 representing the tail risk level.

The theory for static CVaR was developed by B¨auerle and Ott
\cite{BO}, and nested CVaR was studied by
Ruszczy´nski~\cite{Ru10} and Shapiro~\cite{Sh21}. Using the CVaR decomposition theorem by Pflug and Pichler
\cite{PP}, Chow et al.~\cite{CTMP} introduced an RMDP with the risk-augmented states for computing optimal CVaR values by
value iteration. Approximately 10 years later Hau et al.~\cite{Hau23} proved that the result of the value iteration may
not be equal to the optimal value of the static CVaR, and it is a lower bound for the minimal static CVaR.

This paper, which deals with problems with finite state and action sets, introduces a Dynamically augmented
Conditional Value-at-Risk (DCVaR) defined by using the RMDP introduced in Chow et al.~\cite{CTMP},
which we call the Dynamically augmented RMDP (DRMDP). The main result is Algorithm DCVaR for
constructing an optimal policy minimizing the DCVaR value for a given initial tail risk level.

The first group of the results in this paper establishes the relation between the static CVaR values and
the DRMDP. First we observe that there is an optimal nonrandomized policy minimizing static CVaR, and, for a
nonrandomized policy for a Markov Decision Process (MDP), its static CVaR is equal to the worst possible expected outcome if the Decision Maker (DM) plays this policy in
the DRMDP. This observation clarifies the gap described in Hau et al.~\cite{Hau23}
and illustrates time inconsistency
of the static CVaR for MDPs because, in order to implement this worst possible outcome, the second player
in the DRMDP, whom we call Nature, may have to play a policy, which knows future decisions of the
DM. A more natural assumption is that Nature plays its optimal policy, and this assumption leads to the
definition of the DCVaR provided in this paper.  This definition implies that value iterations proposed in Chow et al.~\cite{CTMP} converge to minimal values of DCVaR.

The second group of the results  in this paper deals with formulating Algorithm DCVaR for
constructing optimal DCVaR policies. The proof that this algorithm constructs optimal DCVaR policies is based
on the structure of solutions to a special mass transfer problem describing Nature's optimal decisions.

The optimality equations for the DCVaR is relevant to the optimality equation for the nested CVaR,
and the nested CVaR can be interpreted as the DCVaR with the fixed tail risk level, while optimal DCVaR
policies make decisions for variable tail risk levels.  The tail risk levels for optimal DCVaR policies  depend on the initial tail risk level
and gains and losses incurred up to the current time epoch.  The DM knows only the initial risk level and does not observe risk levels
at later stages, but Algorithm DCVaR implicitly estimates current risk levels.

Classic methods for MDPs deal with risk-neutral objectives. In addition,
there is a significant literature on risk-sensitive optimization initiated by Howard and Matheson~\cite{HM}
dealing mostly with optimization of expectations of utility functions rather than with optimization of
expected total costs or average costs per unit time. Pioneering work by Markowitz~\cite{Markowitz52} on mean-variance
optimization for portfolio selection influenced research on MDPs with mean-variance criteria
\cite{FKL89, Sob}, but
variance of the total costs is not a convenient characteristic to deal with for sequential problems, and it
is not a good risk measure, because variance is affected by losses and gains in the same way. CVaR is an attractive alternative to variance both because of its practical importance and
rich mathematical properties; see \cite{RU00, RU02, RU13, ShDR}. CVaR is a risk measure that has gained popularity in various
applications including finance. It has the ability to safeguard a DM from the tail events
by focusing on an average of the largest losses and by doing so it provides a more comprehensive view for
risk management than threshold-based risk measures such as Value-at-Risk(VaR). CVaR is also known as
expected shortfall in financial literature and has great importance in terms of financial regulation.

We summarize the organization of this paper and its main results. Preliminary information on CVaR and
MDPs is provided in Section~\ref{s2}. In Section~\ref{s3} we prove the existence of a nonrandomized policy minimizing
 CVaR of the total discounted costs. Then we define the DRMDP as an RMDP in which the DM chooses actions, and Nature assigns tail risk levels.  States of the DRMDP are states of the original MDP augmented with tail risk levels.
 Policies for the original MDP are particular policies of the DM in the DRMDP when the DM does not observe past and current tail risk levels, except for the initial one, and we call these policies risk-independent.
 We show that for a nonrandomized risk-independent policy played by the DM,
 CVaR of total discounted costs is equal to the worst expected total discounted costs the DM can incur by playing this policy in the DRMDP. This result illustrates time inconsistency  of the static CVaR and implies that the minimal static
CVaR value for the MDP is equal to the minimal value for the class of nonrandomized risk-independent policies
in the DRMDP. This equality illustrates a positive gap between the optimal static CVaR value and the limit of value iterations. This gap was discovered in Hau et al.~\cite{Hau23}.  Section~\ref{s4} defines a DCVaR, describes
its properties, and introduces a DRMDP1, which is obtained from the DRMDP by modifying one-step costs and
transition probabilities. The payoff and value functions for the DRMDP1 are equal to these functions for the
DRMDP multiplied by  tail risk levels. The sets of optimal policies for the DRMDP and DRMP1 coincide if the
initial tail risk level is positive. The advantage of the DRMDP1 is that its value function is  concave in
the tail risk level.

Section~\ref{s5} introduces the algorithm for constructing a nonrandomized optimal policy minimizing DCVaR,
which was announced in~\cite{DF22}. Sections~\ref{s6}--\ref{s8}
deal with proving that this algorithm constructs an optimal policy. Section~\ref{s6}
describes the properties of the mass transfer problem relevant to optimal solutions for Nature. Section~\ref{s7}
describes the properties of optimal values and properties of the sets of optimal actions. Section~\ref{s8} proves correctness of
the algorithm. Section~\ref{s9} provides extensions to random cost. 
Appendix~\ref{a1} provides results on RMDPs used in this paper. 


\section{Preliminaries: CVaR and MDPs}\label{s2}

Conditional Value at Risk (CVaR) is a  risk measure widely used in engineering and finance. For a  random variable Z defined on a probability space $(\Omega,\mathcal{F},P),$ its CVaR for a tail risk level $\alpha\in [0,1]$ is a conditional expectation of the tail.  One of  several equivalent formal definitions of CVaR is
\begin{equation*}\label{edrcvar} {\rm CVaR}_\alpha(Z): = \frac{1}{\alpha} \int_0^\alpha {\rm VaR}_\beta(Z)d\beta,\qquad \alpha\in (0,1),\end{equation*}
and the Value at Risk is
\[{\rm VaR}_\beta(Z):=\min\{z:F_Z(z)\ge 1-\beta\},\qquad \beta\in (0,1),\]
where $F_Z(z)=P\{Z\le z\}$ is the distribution function of $Z.$  In addition, ${\rm CVaR}_1(Z):= \mathbb{E}[ Z],$ if this expectation exists, and
 ${\rm CVaR}_0(Z):= {\rm ess\, sup}(Z)=\inf\{z\in\R: F_Z(z)=1\},$ where $\min\{\emptyset\}=+\infty$ by definition;  see \cite{RU00,RU02,RU13} and \cite[Chapter 6]{ShDR} for additional details regarding CVaR.

 It is well known that CVaR can be equivalently written in other forms.  For $\alpha\in (0,1]$
\begin{equation*} \label{ecvmin}{\rm CVaR}_\alpha(Z) = \min_{w\in \R}\{w+\frac{1}{\alpha}\mathbb{E}[(Z-w)^{+}]\},\end{equation*}
where $z^+=\max\{z,0\}$ for a number $z.$
{\rm CVaR} can also be represented using its dual representation
\begin{equation*} \label{ecvmax} {\rm CVaR}_\alpha(Z) = \max_{\xi\in\mathcal{U}_{\rm CVaR}(\alpha,P)} \mathbb{E}[\xi Z],\end{equation*}
where $\alpha\in [0,1]$, and $\mathcal{U}_{\rm CVaR}(\alpha,P)$ is the set of random variables $\xi$ on the probability space $(\Omega,\mathcal{F},P)$ defined as $\mathcal{U}_{\rm CVaR}(\alpha,P): = \{\xi\ge 0\, :\,  \alpha\xi\le 1\ {\rm and}\  \mathbb{E}[\xi] = 1\},$ and this set is called the dual CVaR risk envelope.  In this paper, we sometimes write a constant instead of a function equal to this constant, and comparisons between random variables are usually understood $P$-a.s.

Let $\mathcal{G}$ be a sub $\sigma$-algebra of $\mathcal{F},$ that is,  $\mathcal{G}$ is a $\sigma$-algebra on $\Omega$ and $\mathcal{G}\subset \mathcal{F}.$ A $\mathcal{G}$-measurable random variable $\gamma:\Omega\to [0,1]$ is called a random tail risk level.  Then the conditional CVaR given the $\sigma$-algebra $\mathcal{G}$ and a random tail risk level $\gamma$ is defined \cite[Definition 17]{PP} as the $\mathcal{G}$-measurable random variable
\[{\rm CVaR}_\gamma(Z|\mathcal{G}):={\rm ess\, sup}\left\{\E(\xi Z|\mathcal{G}): \xi\in \mathcal{U}_{\rm CVaR}(\gamma,P|\mathcal{G})   \right\}, \]
where the essential supremum is taken with respect to the conditions probability $P|\mathcal{G},$ and $\mathcal{U}_{\rm CVaR}(\gamma,P|\mathcal{G})$ is the set of nonnegative $\mathcal{G}$-measurable random variables $\xi$ on $(\Omega,\mathcal{F},P)$ such that $\E[\xi|\mathcal{G}]=1$ and $\xi\gamma\le 1.$
According to  the CVaR decomposition theorem by   Pflug and Pichler  (\cite{PP}, Lemma 22), 
\begin{equation*}\label{edecom}
{\rm CVaR}_\alpha (Y)=\sup \{ \mathbb{E}[\xi\cdot {\rm CVaR}_{\alpha \xi}(Y|\G)]: \xi\ge 0\ {\rm is}\ \G{\rm -measurable},\ \alpha\xi\le 1,\ {\rm and}\ \mathbb{E}[\xi]=1\}.
\end{equation*}

This paper deals with minimizing CVaR of a random variable $Z_N$ defined as
\begin{equation}\label{eqZdefeff} Z_N:= \sum_{t=0}^N  \beta^t C_t, \end{equation} 
where $(C_t)_{t=0,1,\ldots},$ are random variables, and either $N=1,2,\ldots$ and $\beta\in [0,1],$ or $N=\infty$ and $\beta\in [0,1).$
The constant $\beta$ is called a discount factor.
In particular, this paper deals with the situation when the stochastic sequence $(C_t)_{t=0,1,\ldots}$ can be controlled.  This means that the probability measure $P$ can depend on a policy chosen by a DM, and the goal is to choose the best policy.  The theory of Markov Decision Processes (MDPs) provides a modeling framework for such problems.

An MDP is a model of a natural family of controlled stochastic sequences. It is a tuple $(\X,\A,A(\cdot),c,p),$ where $\X$ is a set of states, $\A$ is a set of actions, and for each state $x\in \X$ there is a nonempty set of available actions $A(x)\subseteq\A,$ $c$ is a real-valued cost function, and $p$ is a transition probability. The time $t=0,1,\ldots$ is discrete,  and, if at some time instance an action $a\in A(x)$ is selected at a state $x\in\X$, then the system moves to the next state $x'\in\X$ according to $x'\sim p(\cdot |x,a)$ and the cost $c(x,a,x')$ is collected, allowing the cost to depend on the next state. In this paper we always assume that $\X:=\{1,\ldots,M\}$ and $\A$ are finite sets. 

Let $\H_t:=(\X\times \A)^t\times\X,$ where $t=0,1,\ldots,$ be the set of histories up to epoch $t=0,1,2,\ldots,$ and $(\X\times\A)^\infty$ is the set of trajectories. This set is endowed with the $\sigma$-algebra $\mathcal{B}((\X\times \A)^\infty),$ which is the countable product of $\sigma$-algebras  $\mathcal{B}(\X\times\A)$  consisting of all subsets of the finite sets $\X\times\A.$

A policy $\pi$ is a sequence of transition probabilities $(\pi_t)_{t=0,1,\ldots}$ from the sets of finite histories $\H_t$ to the sets of actions $\A$ such that
$\pi_t (A(x_t)|h_t)=1$ for all $h_t=x_0,a_0,x_1,a_1,\ldots,x_t\in \H_t.$ Let $\Pi$ be the set of all policies.   According to the Ionescu Tulcea theorem, an initial state $x\in\X$ and  a policy $\pi\in\Pi$ define  a probability $P_x^\pi$ on the space of trajectories $(\X\times\A)^\infty, $ where
$\pi_t$ are transition probabilities from $\H_t$ to $\A,$ ana $p$ is a transition probability from $\X\times\A$ to $\X,$ which can be viewed as transition probabilities from $\H_t\times \A$ to $\X,$ $t=0,1,\ldots.$ Expectations with respect to probabilities $P_x^\pi$ are denoted by  $\E_x^\pi.$
If each probability $\pi(d_t|h_t)$ is concentrated at one point, the policy $\pi$ is called nonrandomized.  A nonrandomized policy is defined by a sequence of mappings $\phi_t:\H_t\to\A$ such that $\phi_t(x_0,a_0,\ldots x_{t-1},a_{t-1},x_t)\in A(x_t)$ for all $t=0,1,\ldots$ and for all $x_t\in\X.$ A nonrandomized policy is called Markov if each time the choice of a decision depends only on the current time and state.   A Markov policy is defined by a sequence of mappings $\phi_t:\X\to\A$ such that $\phi_t(x)\in A(x)$ for all $t=0,1,\ldots$ and for all $x\in\X.$  A Markov policy is called deterministic if decisions do not depend on the time parameter.  A deterministic policy is defined by a function $\phi:\X\to \A$ such that $\phi(x)\in A(x)$ for all $x\in\X.$
Let $\beta\in [0,1]$ be a constant discount factor. 

 An objective criterion $v(x,\pi)$ is a real-valued function on $\X\times\Pi.$  In general, it can take infinite values, but in this paper we consider only real-valued objective criteria. The function  $v(x):=\inf_{\pi\in \Pi} v(x,\pi),$ where $x\in\X,$ is called the value function.   A policy $\pi$ is optimal for the initial stated $x\in\X$ if $v(x,\pi)=v(x).$  A policy $\pi$ is called optimal if it is optimal  for all initial states $x\in\X.$  In the literature these definitions are usually applied to risk-neutral objectives, but they can be applied for arbitrary  objective functions \cite{F82} with values in $\R\cup\{\pm\}.$

 For a finite horizon $N=1,2,\ldots$ the total discounted cost $Z_N$ mentioned in \eqref{eqZdefeff}    is the random variable
\begin{equation}\label{eqCN}
Z_N:= \sum_{t=0}^{N-1} \beta^t c(x_t,a_t,x_{t+1})+\beta^Nv_0(x_N),
\end{equation}
where $v_0:\X\to\R$ is the terminal cost, and for the
 infinite-horizon $N=\infty$ 
\begin{equation}\label{eqC}
Z_\infty = \sum_{t=0}^\infty \beta^t c(x_t,a_t,x_{t+1}),
\end{equation}
where $\beta\in [0,1),$ if $N=\infty,$ and $\beta\in [0,1]$ if $N<\infty$.  In fact, for $N<\infty,$ the results of this paper holds for nonnegative real-valued discount factors $\beta.$

The risk-neutral approach is to consider the expected total discounted cost
$\E_x^\pi [Z_N].$
For a tail risk level $\alpha\in [0,1]$ and  a finite horizon $N=1,2,\ldots$ or an infinite horizon $N=\infty,$ let us consider the objective function  $ {\rm CVaR}_\alpha(Z_N;P_x^\pi),$ where $x\in\X,$ $\pi\in\Pi,$ , and $Z_N$ is the random variable defined on the probability space $(\Omega,\mathcal{F},P)=((\X\times\A)^\infty, \mathcal{B}((\X\times\A)^\infty),P_x^\pi)$ in formulae \eqref{eqCN} or \eqref{eqC}. The value ${\rm CVaR}_\alpha(Z_N;P_x^\pi)$ is sometimes called a static CVaR.  Computing optimal policies minimizing the static CVaR is a hard problem  \cite[Remark 5.1]{BO}. 

Chow et al.~\cite{CTMP} changed the original problem formulation with the state space $\X$ to the formulation when the states of the problem consist of state-risk pairs $(x,y),$ where $x\in\X$ is the original state and $y\in\ [0,1]$ is the tail risk level.  In other words, the state space $\X$ is augmented in \cite{CTMP} with the tail risk level, and the new state space is  ${\bf X}:=\X\times [0,1].$  Motivated by the Pflug-Pichler CVaR decomposition theorem \cite{PP}, Chow et al.~\cite{CTMP} introduced the Robust MDP (RMDP) with the state space ${\bf X},$ which is formally defined in Section~\ref{s3} and called the Dynamically augmented RMDP (DRMDP). Infinite horizon was considered in \cite{CTMP}, and a finite horizon version was considered in \cite{LZB}, where it was assumed that the DM knows tail risk levels at all steps.  Hau et al.~\cite{Hau23} showed that the value of DRMDP is a lower bound of the minimal CVaR values for the original MDP.  In Section~\ref{s4} we introduce the Dynamically augmented CVaR (DCVaR) objective function, which is a lower bound of CVaR, and the value of DRMDP is equal to the minimal value of DCVaR.

In the conclusion of this section, we recall two concepts for MDPs relevant to more general objective criteria than expected total costs: MDPs with general objective functions and abstract dynamic programming. According to \cite{F82}, a  general objective function is a function $v:\X\times\Pi\to \R\cup \{\pm\infty\}$ such that, if $P_x^\pi=P_x^\sigma$ for $x\in\X$ and $\pi,\sigma\in\Pi,$ then $v(x,\pi)=v(x,\sigma).$ This definition is consistent with the notion of a law-invariant risk measure \cite{PP, ShDR}.

Abstract dynamic programming~\cite{Bert} is a classic approach based on representing the optimality equation in the general form $v(x)=\inf_{a\in A(x)}H(x,a,v).$ This approach is relevant to the notion of the nested CVaR~\cite{ShDR}.

\section{Static CVaR and Robust MDPs}\label{s3}
The main results of this section are: (i)  the existence of nonrandomized policies minimizing static CVaR (Theorem~\ref{thexoptcvar}), and (ii) the minimal value of  static CVaR is the minimal  guaranteed expected total discounted payoff for the DM in the RMDP introduced in Chow et al.~\cite{CTMP}, if the DM plays policies $\pi\in\Pi$ from the MDP; see  Corollary~\ref{cor34} below. This corollary provides   a game-theoretic interpretation of the gap between the optimal value of CVaR and optimal value of this RMDP discovered in Hau et al.~\cite{Hau23}.


Let us consider an MDP $(\X,\A,A(\cdot),c,p)$ with finite state and action spaces $\X$ and $\A$ respectively. The goal is to find a policy $\pi$ minimizing ${\rm CVaR}_\alpha(Z_N;P_x^\pi)$ for all $x\in\X$ and for a given risk level $\alpha\in [0,1].$
Let \[{\rm CVaR}_\alpha(Z_N;x):=\inf_{\pi\in\Pi} {\rm CVaR}_\alpha(Z_N;P_x^\pi)\] be the optimal value of the static CVaR.  The following theorem states the existence of nonrandomized optimal policies minimizing static CVaRs and convergence of finite-horizon optimal values of static CVaRs to the optimal infinite-horizon value. We recall that $\alpha\in [0,1]$ is fixed. We also recall that throughout the entire paper the discount factor $\beta\in [0,1],$ if $N=1,2,\ldots,$ and $\beta\in [0,1)$ if $N=\infty.$
\begin{theorem}\label{thexoptcvar}
For every $N=1,2,\ldots$ or for $N=\infty,$ for the static CVaR optimization problem there exists a nonrandomized optimal policy $\phi\in\Pi$  for which ${\rm CVaR}_\alpha(Z_N;P_x^\phi)={\rm CVaR}_\alpha(Z_N;x)$ for all $x\in\X.$  In addition,
\begin{equation*}\label{eqlimefgr}
{\rm CVaR}_\alpha(Z_\infty;x)=\lim_{N\to\infty} {\rm CVaR}_\alpha(Z_N;x),\qquad x\in\X.
\end{equation*}
\end{theorem}
\begin{proof}
We observe that it is sufficient to proof that for every $x\in\X$ there exists a nonrandomized optimal policy $\phi^x$ for the initial state $x.$ Indeed, let $\Delta^N$ be the set of all nonrandomized $N$-horizon policies.   If  $\phi^x\in \Delta^N$ is optimal for an initial state $x\in\X,$ then the nonrandomized policy $\phi$ is optimal, where
$
\phi_t(h_t):= \phi^{x_0}(h_t), \qquad t<N,\ h_t=x_0,a_0,\ldots,x_t.
$
Let us fix $x\in\X$ and prove the existence of $\phi^x.$ Let us start with $N<\infty.$ In this case the set $\Delta^N$ is finite.  Therefore, the strategic measure for an arbitrary policy $\sigma$ can be presented as a convex combination of strategic measures for nonrandomized policies \cite[Theorem 1]{F82a}: there are numbers $\lambda(\pi)\ge 0$ such that $\sum_{\pi\in\Delta^N}\lambda(\pi)=1$ and
$P_x^\sigma=\sum_{\pi\in \Delta^N} \lambda(\pi) P_x^\pi.$ Since taking CVaR of a mixture of distributions is a concave operation \cite[Proposition 3.2]{PU},
\[{\rm CVaR}_\alpha (Z_N;P_x^\sigma)= {\rm CVaR}_\alpha (Z_N;\sum_{\pi\in \Delta^N}\lambda(\pi) P_x^\pi)\ge \sum_{\pi\in \Delta^N}\lambda(\pi) {\rm CVaR}_\alpha (Z_N;P_x^\pi)\ge {\rm CVaR}_\alpha (Z_N;P_x^{\phi^x}),\]
where $\phi^x\in\Delta^N$ such that ${\rm CVaR}_\alpha (Z_N;P_x^{\phi^x})=\min_{\pi\in\Delta^N}  {\rm CVaR}_\alpha (Z_N;P_x^\pi). $ Thus, for an arbitrary policy $\sigma$ we have that ${\rm CVaR}_\alpha (Z_N;P_x^{\phi^x})\le {\rm CVaR}_\alpha (Z_N;P_x^\sigma),$ which means that that a deterministic policy $\phi^x$ is optimal for the $N$-horizon problem with the initial state $x.$

Let us consider the infinite-horizon problem. In the rest of this proof, let us set $v_0(x):=0$ for all $x\in\X.$ If we add a constant $d$ to the cost function $c,$ then the all CVaR values for an $N$-horizon problem will be changed by the constant $d(1-\beta^N)/(1-\beta).$  Thus, this addition preserves sets of optimal policies.  Since the cost function $c$ is bounded, without loss of generality we can assume that the function $c$ is nonnegative.  Let us assume that $c(x,a,x')\ge 0,$  and let $K>0$ is an upper bound of the cost function, that is,  $0\le c(x,a,x')\le K $ for all $x,x'\in\X$ and $a\in A(x).$

Let $N<\infty$ and $\tilde{N}>N.$   Then for $Z_{\tilde{N}}\ge Z_N$ and $Z_{\tilde{N}}-Z_N\le Z_\infty-Z_N\le \epsilon_N\to 0$ as $N\to\infty$ for  all trajectories, where $\epsilon_N= K\beta^N/(1-\beta).$  Therefore, for every policy $\pi\in\Pi$
\begin{equation}\label{eqcomptwopol}
{\rm CVaR}_\alpha(Z_N;P_x^\pi)\le {\rm CVaR}_\alpha(Z_{\tilde{N}};P_x^\pi)\le {\rm CVaR}_\alpha(Z_N;P_x^\pi)+\epsilon_N,
\end{equation}
which implies
\[
{\rm CVaR}_\alpha(Z_N;x)\le {\rm CVaR}_\alpha(Z_{\tilde{N}};x)\le {\rm CVaR}_\alpha(Z_N;x)+\epsilon_N,
\]
and

\begin{equation}\label{eqlimnvef}
 {\rm CVaR}_\alpha(Z_N;P_x^\pi)\uparrow {\rm CVaR}_\alpha(Z_\infty;P_x^\pi), \quad  {\rm CVaR}_\alpha(Z_N;x)\uparrow {\rm CVaR}_\alpha(Z_\infty;x)\  {\rm as}\  N\to\infty.
 \end{equation}

Let us construct a nonrandomized optimal infinite-horizon policy $\phi\in\Delta^\infty.$  Let $\phi^N\in\Delta^N,$ $N=1,2,\ldots,$ be nonrandomized optimal $N$-horizon policies. Since the sets $\X$ and $A(x),$ where $x\in\X,$ are finite, it is possible to choose a sequence $\{N^1_i\}_{i=1}^\infty$ of integers such that  $\phi_0^{N^1_i}(x_0)=\phi_0^{N^1_j}(x_0)$ for all $x_0\in\X$ and for all $i,j=1,2,\ldots.$ We define
$\phi_0(x_0):=\phi_0^{N^1_i}(x_0)$ for all $x_0\in\X.$

Now let for some $n=1,2,\ldots$ there is an increasing sequence of natural numbers $\{N^n_i\}_{i=1}^\infty$ such that for each history $h_t=x_0,a_0,x_1,\ldots x_t,$ $t=0,1,\ldots,n-1,$ policies $\phi^{N^n_i}$ make the same decisions $\phi_t^{N^n_i}(h_t)$ for all $i=1,2,\ldots.$   We select a subsequence  $\{N^{n+1}_i\}_{i=1}^\infty$ of the sequence $\{N^n_i\}_{i=1}^\infty$ such that for each history $h_n=x_0,a_0,x_1,\ldots,x_n$ all policies $\phi^{N^{n+1}_i}$ make the  decision $\phi_n^{N^{n+1}_i}(h_n)$ at the step $n,$ which is the same for all $i.$ 
By repeating this procedure and taking $n\to\infty,$ we construct a nonrandomized policy  $\phi.$

Let us show that $\phi$ is indeed an optimal policy.  Let us consider a time horizon $n.$  Then there exists an integer $\tilde{N}=N^{n+1}_i>n$ such that the policy  $\phi[\tilde{N}]$ is optimal for the horizon $\tilde{N}$ and coincides at the first $n$ steps with the policy $\phi$. Therefore, for this policy $\phi$, \eqref{eqcomptwopol} an be rewritten as
\begin{equation*}
{\rm CVaR}_\alpha(Z_n;P_x^\phi)\le {\rm CVaR}_\alpha(Z_{\tilde{N}};P_x^{\phi[\tilde{N}]})\le {\rm CVaR}_\alpha(Z_n;P_x^\phi)+\epsilon_n.
\end{equation*}
Since the policy $\phi[\tilde{N}]$ is optimal for the horizon $\tilde{N},$
\begin{equation*}
{\rm CVaR}_\alpha(Z_n;P_x^\phi)\le {\rm CVaR}_\alpha(Z_{\tilde{N}};x)\le {\rm CVaR}_\alpha(Z_n;P_x^\phi)+\epsilon_n.
\end{equation*}
Let $n\to\infty.$  Then $\epsilon_n\to 0.$ Since $\tilde{N}\ge n,$ in view of \eqref{eqlimnvef},
\[
{\rm CVaR}_\alpha(Z_\infty;P_x^\phi)= {\rm CVaR}_\alpha(Z_\infty;x), \qquad \alpha\in [0,1],\ x\in\X.
\]
\end{proof}

We recall that $\X=\{1,\ldots,M\}$ and denote by $\R^M$ the $M$-dimensional Euclidian space. Let us consider the RMDP introduced in Chow et al.~\cite{CTMP}, which we call the Dynamically augmented RMDP (DRMDP).  The DRMDP is defined by a tuple
$({\bf X},\A,\B,A(\cdot),B(\cdot,\cdot,\cdot),c,q),$ where the state space is ${\bf X}:=\X\times [0,1],$ action space is $\A,$ uncertainty space is $B:=\R^M$ with $M$ being the number of states in $\X,$ action sets for the DM at states $(x,y)\in\bf{X}$ are $ A(x,y):=A(x),$   uncertainty sets for Nature are
\begin{equation*}\label{UCVaR1}
B (x,y,a):=\mathcal{U}(x,y,a)\cap \{b\in \R^M_+:\, b_{x'}=0 \ {\rm if}\ p(x'|x,a)=0 \},
\end{equation*}
where, $\R^M_+:=\{b\in\R^M\, : b_{x'}\ge 0, x'\in\X\},$ and, for $x\in\X,$ $y\in [0,1],$ $a\in A(x),$
\begin{equation*}\label{UCVaR1sup}
\mathcal{U}(x,y,a):=\{b\in {\R}^M_+ :   yb_{x'}\le 1, x'=1,2,\ldots,M, \sum_{x'=1}^M b_{x'}p(x'|x,a)=1\},
\end{equation*}
  one-step costs $c((x,y),a,b,(x',y')):=c(x,a,x')$  for $(x,y),(x',y')\in {\bf X},$ $a\in A(x),$ $b\in B(x,y,a,),$  and  transition probabilities
$q(x',D|x,y,a,b):=b_{x'}p(x'|x,a)\delta_{yb_{x'}}(D)$ for $x,x'\in\X,\ y\in [0,1],$ $D\in\mathcal{B}([0,1]),$ $a\in A(x),$ and  $b\in \mathcal{U}(x,y,a),$
where $\delta_z(\cdot)$ is the Dirac measure on the interval $[0,1]$ concentrated at the point $z\in [0,1].$  A DRMDP is a particular case of the RMDP model described in Appendix~\ref{a1}, when  the sets ${\bf X},$ $\B,$ set-valued mapping $B,$ one-step cost function $c,$ and transition probability $q$ have specific forms.  In particular, states for DRMDPs are denoted by $(x,y)\in{\bf X}.$ The sets of policies for the DM and Nature are denoted by $\Pi^\A$  and $\Pi^\B$ respectively. We shall use notation  $v_N(x,y,\pi^\A,\pi^\B):=v_N((x,y),\pi^\A,\pi^\B),$ where $(x,y)\in{\bf X},$ $\pi^\A\in \Pi^\A,$ and $\pi^\B\in \Pi^\B,$  for the expected total discounted payoffs over the horizon $N$ with the terminal payoff $v_0(x,y),$ which is continuous in $y\in [0,1].$  As explained after formulae \eqref{eqsecondopa} and \eqref{equnbanch}, the optimal values $v_N:\X\times [0,1]\to \R$ are continuous functions.  Since $\X$ is a finite set, this means that the functions $v_N(x,y)$ are continuous in $y\in [0,1].$

 We observe that $\Pi\subset \Pi^\A,$ that is, the  set of all policies $\Pi$ for the original MDP is the subset of the set of policies for the DM in the DRMDP.  We call policies from $\Pi$ risk-independent. The fundamental difference between a policy $\pi\in \Pi$ and a policy from $\Pi^\A$ is that the policy $\pi$ does not have information about the current and past tail risk levels $y_t$ except the initial tail risk level $y_0=\alpha,$ which we consider to be fixed. 

The following theorem states that, for each risk-independent policy for the DM, Nature can play a policy maximizing DM's losses.  The proofs of Theorems~\ref{thoptalbhd} and \ref{thmA2} are based on the well-known fact that, if in a two-player stochastic game a policy of one player is fixed, then another player deals with an MDP.  In general, the states of this MDP are histories in the game.  However, the state space of this MDP can be simplified, if the fixed policy belongs to a special class.  For example, if the fixed policy is stationary, then another player deals with an MDP whose states are the states of the stochastic game.  There are differences between the proofs of Theorems~\ref{thoptalbhd} and \ref{thmA2}. The proof of Theorem~\ref{thmA2} follows from the validity of sufficient conditions for the existence of optimal policies for MDPs with setwise continuous transition probabilities \cite{Sch75,FKZ}, which takes place because all the sets $A(x)$ are finite. The proof of Theorem~\ref{thoptalbhd} follows from the validity of sufficient conditions for the existence of optimal policies for MDPs with weakly continuous transition probabilities \cite{Sch75,FK20}, which takes place because the set $\X$ and all the sets $A(x)$ are finite and because of weak continuity of the transition probability $q.$ The remaining statements of this section deal with the DRMDP.

\begin{theorem}\label{thoptalbhd} For each  risk-independent policy $\pi\in\Pi$ for the DM and for each $N=1,2,\ldots$ or $ N=\infty,$  there exists a nonrandomized policy
$\phi^\B\in\Pi^\B,$ which depends on $N$ and $\alpha,$ such that
\[
v_N(x,\alpha,\pi,\phi^\B)=\max_{\pi^\B\in\Pi^\B} v_N(x,\alpha,\pi,\pi^\B)\quad {\rm for\ all\ } x\in\X,\ \alpha\in [0,1].
\]
\end{theorem}
\begin{proof}
If the DM plays a risk-independent policy $\pi\in\Pi,$ Nature deals with an MDP with states $\tilde{x}_t:=x_0,a_0,x_1,\ldots,x_t,y_t,a_t\in  (\X\times \A)^t\times\X\times [0,1]\times\A,$ $t=0,1,\ldots.$ The set of available actions at each state $\tilde{x}_t$ is $B(\tilde{x}_t):= B(x_t,y_t,a_t),$ where $x_t\in\X, $ $y_t\in [0,1],$ and $a_t\in A(x_t).$ The one-step reward for Nature is $c(\tilde{x}_t, b_t, \tilde{x}_{t+1}):=c(x_t,a_t,x_{t+1}),$ and this function is bounded and continuous, where $\tilde{x}_{t+1}:= x_0,a_0,x_1,\ldots,x_t,a_t, x_{t+1},y_{t+1},a_{t+1}.$  If an action $b$ is chosen at state  $\tilde{x}_t,$ then the next state is $\tilde{x}_{t+1}$  with the probability
$\tilde{p}(\tilde{x}_{t+1}|\tilde{x}_t,b):=b_{x_{t+1}}p(x_{t+1}|x_t,a_t)\pi_{t+1}(a_{t+1}|x_0,a_0,\ldots,x_{t+1}),$  where the following conditions hold: (a) $\tilde{x}_{t+1}=x_0,a_0,x_1,\ldots,x_t,a_t,x_{t+1},y_{t+1},a_{t+1},$ (b) $x_{t+1}\in\X,$ (c) $y_{t+1}=y_tb_{x_{t+1}},$ (d) $a_{t+1}\in A(x_{t+1}).$ Other transitions are impossible since
$\sum_{x'\in\X}b_{x'}p(x'|x,a)=1.$    The transition probability $\tilde{p}(d\tilde{x}_{t+1}|\tilde{x}_t,b)$ is weakly continuous in $(\tilde{x}_t,b)$ because $b\mapsto b_{x_{t+1}}$ is a bounded continuous function and $p(x_{t+1}|x_t,a_t)$ and $\pi_{t+1}(a_{t+1}|x_0,a_0,\ldots,x_{t+1})$ are distributions on finite sets depending on the finite numbers of conditions.  We consider the expected total discounted rewards for this problem with the discount factor $\beta.$  The initial state is $\tilde{x}_0=(x_0,y_0,a_0),$ where $y_0=\alpha,$ which is the initial tail risk level, and $a_0\sim \pi_0(x_0|x_0).$ This MDP satisfies the weakly continuous conditions, which implies the existence of optimal Markov policies for finite-horizon problems and optimal deterministic policies for infinite-horizon problems \cite{Sch75} or \cite[Theorem 2]{FKZ}. These optimal policies define the policy $\phi^\B$ for Nature whose existence is stated in the theorem.
\end{proof}


The following theorem establishes the relations between the DRMDP   and  CVaR of the total discounted reward defined in \eqref{eqZdefeff}. In particular, formula \eqref{eqth33} illustrates time inconsistency of static CVaR because at some time instance the policy of Nature maximizing the right-hand side of \eqref{eqth33} may depend on future decisions of the DM playing a possibly nonstationary policy $\phi.$

\begin{theorem}\label{thoptalbhd1} For each nonrandomized risk-independent policy $\phi\in\Pi$ for the DM, for each initial tail risk level $\alpha\in [0,1]$, and for each $N=1,2,\ldots$ or $N=\infty,$  
\begin{equation}\label{eqth33}
{\rm CVaR}_\alpha(Z_N;P_x^\phi)=\max_{\pi^\B\in\Pi^\B} v_N(x,\alpha,\phi,\pi^\B),\qquad x\in\X.
\end{equation}
\end{theorem}
\begin{proof}
Let us consider the MDP defined at the beginning on the proof of Theorem~\ref{thoptalbhd} for $\pi=\phi.$  We  fix $N<\infty$ and compute an optimal policy for Nature. If at time $t=N-1$ Nature is at state $\tilde{x}_{N-1}$ of this MDP, then the optimal 1-step value for Nature is
\begin{equation*}\label{eqv1forBN}
\hat{V}_1(\tilde{x}_{N-1})=\max_{b\in B(x_{N-1},y_{N-1},a_{N-1})}\sum_{x_N\in\X}[c(x_{N-1},a_{N-1},x_N)+\beta \hat{V}_0(x_N,y_{N-1}b_{x_N})]b_{x_N}p(x_N|x_{N-1},a_{N-1}),
\end{equation*}
 where $\hat{V}_0(x,y)=v_0(x,y),$ $x\in\X,$ $y\in [0,1],$ there is a nonrandomized optimal 1-step Markov policy $\phi_{N-1}^\B(x_{N-1},y_{N-1},a_{N-1})$ \cite[Theorem 2]{FKZ}, and, in view of Berge's maximum theorem, the function $\hat{V}_1(\tilde{x}_{N-1})$ is continuous. This means that this function is continuous  in $y_{N-1}$ since all other variables in $\tilde{x}_{N-1}$ take a finite number of values.  Nature chooses the optimal action $b_{N-1}=\phi_{N-1}^\B(x_{N-1},y_{N-1},a_{N-1})$ without using any knowledge of the previous states and actions.  This is true for $t=N-1$ because the optimal value of $\hat{V}_1$  does not depend on the distributions of future actions and states following the terminal state $(x_N,y_N).$

Let us make the induction assumption that for some integer  $t=N-2,N-3,\ldots,1,0$ there is a Markov policy $(\phi_{t+1}^\B,\phi_{t+2}^\B,\ldots,\phi_{N-1}^\B)$ such that this policy is optimal for the horizon $(N-t-1),$ and the value function $\hat{V}_{N-t-1}(\tilde{x}_{t+1})$ is continuous, which means that it is continuous in $y_{t+1}\in [0,1].$  Then

\begin{equation}\label{eqvtforBN}
\hat{V}_{N-t}(\tilde{x}_t)=\max_{b\in B(x_t,y_t,a_t)}\sum_{x_{t+1}\in\X}[c(x_t,a_t,x_{t+1})+\beta \hat{V}_{N-t-1}(\tilde{x}_{t+1})]b_{x_{t+1}}p(x_{t+1}|x_t,a_t),
\end{equation}
for  $\tilde{x}_{t+1}= \tilde{x}_t,a_t,x_{t+1},y_{t+1},a_{t+1},$ where $y_{t+1}=y_tb_{x_{t+1}}$ and $a_{t+1}= \phi_{t+1}(x_0,a_0,\ldots,x_{t+1}).$   In view of Berge's maximum theorem, the function  $\hat{V}_{N-t}(\tilde{x}_{t})$ is continuous in $y_t,$ and there is a measurable mapping $ \tilde{x}_t\mapsto \phi_{N-t}^\B( \tilde{x}_t)$ such that the maximum in
\eqref{eqvtforBN} is achieved at $b= \phi_{N-t}^\B( \tilde{x}_t);$ \cite[Theorem 2]{FKZ} or \cite[Theorem 3.3]{FKZ1}. Thus, the Markov policy $(\phi_t^\B,\phi_{t+1}^\B,\ldots,\phi_{N-1}^\B)$ is optimal for the horizon $(N-t).$  By induction, the Markov policy $(\phi_0^\B,\phi_1^\B,\ldots,\phi_{N-1}^\B)$ is optimal for the horizon $N.$  This Markov policy is defined for  Nature in the special MDP.  It corresponds to the nonrandomized policy for Nature, which with a small abuse of notations we denote by  $\phi^\B,$
\[
\phi_t^\B(x_0,y_0,a_0,x_1,y_1,a_1,\ldots,x_t,y_t,a_t):=\phi_t^\B(\tilde{x}_t),\qquad t=0,1,\ldots,N-1.
\]


If the DM plays the nonrandomized risk-independent policy $\phi\in\Pi,$ and Nature plays the nonrandomized policy $\phi^\B,$ then, according to Theorem~\ref{thoptalbhd},
\begin{equation}\label{eqoprespBN}
\hat{V}_N(x,\alpha,\phi_0(x))=v_N(x,\alpha,\phi,\phi^\B)=\sup_{\pi^\B\in\Pi^\B} v_N(x,\alpha,\phi,\pi^\B), \qquad x\in\X, \alpha\in [0,1].
\end{equation}

 A minor change in the described construction leads to defining an optimal nonrandomized policy for Nature, if the DM plays an arbitrary risk-independent policy $\pi\in\Pi,$ but we do not need the explicit form for the policy $\phi^\B$ for Nature if the DM has chosen  a possibly randomized risk-independent policy $\pi.$

In view of the CVaR decomposition theorem~\cite{PP} and the definition of sets $B(x,y,a),$ for $x\in\X$
\[{\rm CVaR}_{y_{N-1}}(Z_N;P_x^\phi|\tilde{x}_{N-1})=\sum_{t=1}^{N-1} \beta^{t-1}c(x_{t-1},a_{t-1},x_t)+\beta^{N-1}\hat{V}_1(\tilde{x}_{N-1}),\]
for $t=1,\ldots,N-2,$
\[
{\rm CVaR}_{y_t}(Z_N;P_x^\phi|\tilde{x}_t)=\sum_{s=1}^t \beta^{s-1}c(x_{s-1},a_{s-1},x_s)+\beta^t\hat{V}_{N-t}(\tilde{x}_t),
\]
and, for $t=0$ and $\alpha\in [0,1],$ we have
$
{\rm CVaR}_\alpha(Z_N;P_x^\phi)=\hat{V}_N(x,\alpha,\phi_0(x))
.$ The last equality and \eqref{eqoprespBN} imply that the theorem is proved for $N<\infty.$

For $N=\infty,$

\[
{\rm CVaR}_\alpha(Z_\infty;P_x^\phi)=\lim_{N\to\infty}{\rm CVaR}_\alpha(Z_N;P^\phi_x)=\lim_{N\to\infty}\sup_{\pi^\B\in\Pi^\B} v_N(x,\alpha,\phi,\pi^\B)\] \[=\sup_{\pi^\B\in\Pi^\B} v_\infty(x,\alpha,\phi,\pi^\B),
\]
where the first equality follows from unform $P_x^\phi$-a.s. convergence $Z_N\to\Z_\infty,$ the second equality follows from the proved part of this theorem for $N<\infty,$ and the last one follows from convergence of value iterations for the MDP with weakly continuous transition probabilities \cite[Theorem 2]{FKZ} introduced in the proof of Theorem~\ref{thoptalbhd} when $\pi=\phi.$
\end{proof}

The following corollary from Theorems~\ref{thexoptcvar} and \ref{thoptalbhd1} characterizes the optimal value of static CVaR in terms of the DRMDP. Recall that $\Pi$ is the set of policies for the original MDP.
\begin{corollary}\label{cor34} For every $N=1,2,\ldots$ or $N=\infty,$ every $x\in\X$, and every $\alpha\in [0,1],$
\[ {\rm CVaR}_\alpha (Z_N;x)=\min_{\phi^\A\in\Pi_{NR}} \max_{\pi^\B\in\Pi^\B} v_N(x,\alpha,\phi^\A,\pi^\B),\]
where $\Pi_{NR}$ is the set of nonrandomized risk-independent policies for the DM, and $\Pi_{NR}\subset\Pi\subset\Pi^\A.$
\end{corollary}

We recall that the DRMDP satisfies assumptions of an RMDP considered in Appendix~\ref{a1}.  Thus, the value functions $v_N(x,\alpha)$ exist and are continuous in $\alpha\in [0,1]$ for finite and infinite time horizons $N$, and $v_N(x,\alpha)\to v_\infty(x,\alpha)$ as $N\to\infty$ uniformly in $x$ and $\alpha.$ Theorem~\ref{thmA1}  implies the existence of optimal policies and characterizes some classes of optimal policies.  In particular, for each player there exist  optimal nonrandomized Markov policies, and for $N=\infty$ each player has a deterministic optimal policy.

Let us consider the gap between the optimal value of CVaR and the value of the RMDP.  For $N=0,1,\ldots$ or $N=\infty,$ for $x\in\X,$ and for $\alpha\in [0,1],$
\begin{eqnarray*}\Delta_N(x,\alpha):&=&{\rm CVaR}_\alpha(Z_N;x)-v_N(x,\alpha)\\ =\min_{\phi^\A\in\Pi_{NR}} \max_{\pi^\B\in\Pi^\B} v_N(x,\alpha,\phi,\pi^\B)&-&\min_{\pi^\A\in\Pi^\A} \sup_{\pi^\B\in\Pi^\B} v_N(x,\alpha,\pi^\A,\pi^\B)\ge 0,
\end{eqnarray*}
where the second equality follows from Corollary~\ref{cor34} and Theorem~\ref{thmA3}, and the inequality follows from $\Pi_{NR}\subset\Pi\subset \Pi^\A.$ It is shown in \cite{Hau23} that this gap can be positive for $N>1;$ see also \cite{GD}. 

If $\beta=0,$ then we deal with a one-step problem, which can be solved easily by performing one step of value iterations  for the DRMDP with $v_0\equiv 0.$ In this case, the gap is equal to 0.  Therefore, in the rest of this paper, except Appendix~\ref{a1}, we consider $\beta>0.$

\section{Dynamically Augmented CVaR}\label{s4}
Let us consider the DRMDP with $\beta\in (0,1],$ if $N<\infty,$ and $\beta\in (0,1)$ if $N=\infty.$  Theorem~\ref{thoptalbhd1} demonstrates time inconsistency of the static CVaR for MDPs because, for a policy $\pi^\B\in\Pi^\B$ on which CVaR is achieved, current decisions of Nature may depend on future decisions of the DM.  The following definition addresses this problem. Optimal decisions of Nature do not depend on future actions of the DM. Therefore, it is natural for Nature to make optimal decisions.

\begin{definition}
	For a policy $\pi^\A\in\Pi^\A,$ initial state $x\in\X,$  tail risk level $\alpha\in [0,1],$ and time horizon $N=1,2,\ldots$ or $N=\infty,$ the Dynamically augmented CVaR  (DCVaR) is
	\[
	{\rm DCVaR}_\alpha(Z_N;x,\pi^\A):=\sup_{\pi^\B\in\Pi^\B_*(N) }v_N(x,\alpha,\pi^\A,\pi^\B),
	\]
where $	\Pi^\B_*(N)$ is the set of optimal policies for Nature for the $N$-horizon problem.
\end{definition}

  Since $\Pi\subset\Pi^\A,$ this definitions applies to policies $\pi\in \Pi.$ According to Theorem~\ref{thexoptcvar}, for each tail risk level $\alpha\in [0,1]$ there is an optimal nonrandomized policy $\phi\in\Pi$  minimizing the CVaR for the MDP. If the DM plays a nonrandomized risk-independent policy $\phi\in\Pi,$ then, in view of  $\Pi^\B_*(N)\subset \Pi^\B$ and Theorem~\ref{thoptalbhd1}, 
  \[{\rm DCVaR}_\alpha(Z_N;x,\phi)\le {\rm CVaR}_\alpha(Z_N;P_x^{\phi}). \] Thus, in addition to being a more natural objective function than CVaR, DCVaR can be used for establishing performance guarantees for CVaR.    In addition, for the DRMDP
 \begin{equation}\label{eqvalDVARDVAR} v_N(x,\alpha)=\min_{\pi^\A\in\Pi^\A} {\rm DCVaR}_\alpha(Z_N;x,\pi^\A).\end{equation}
 This is true because
 for each $\pi_*^\B\in\Pi_*^\B (N)$
 \begin{equation}\label{eq4.2}
 v_N(x,\alpha)=\min_{\pi^\A\in\Pi^\A}v_N(x,\alpha,\pi^\A,\pi_*^\B)
 \le \min_{\pi^\A\in\Pi^\A} {\rm DCVaR}_\alpha(Z_N;x,\pi^\A)\le v_N(x,\alpha),
 \end{equation}
where the equality follows from the existence and definition of an optimal policy for the DM in the DRMDP, and the inequalities follow from the definition of the DCVaR and from the first equality in Theorem~\ref{thmA3}. Thus all inequalities in \eqref{eq4.2} hold in the form of equalities.  

According to \eqref{eqminimaxgejn}, for $N=1,2,\ldots$ or $N=\infty$ and for $x\in\X,$ $y\in [0,1]$
 \begin{equation}\label{eqVRMDPCV}
v_N(x,y)=\min_{a\in A(x)}\max_{b\in B(x,y,a)} \sum_{x'\in\X}[ c(x,a,x')+\beta  v_{N-1}(x',yb_{x'})]b_{x'}p(x'|x,a).
\end{equation}
We recall that the bounded continuous  function $v_0$ is given,  the functions $v_N(x,y)$ are bounded and continuous in $y,$ and $v_N(x,y)\to v_\infty(x,y)$ uniformly in $y$ as $N\to\infty;$ see \eqref{equnbanch}. Values $v_N(x,y)$ can be computed by value iterations, and the interval $[0,1]$ can be easily discretized for computations.  In addition $v_\infty$ is the unique bounded solution of \eqref{eqVRMDPCV} with $N=\infty.$

Here we would like to mention that time-inconsistency of the static CVaR for MDPs is well-known, and another objective, called the nested CVaR, is broadly used \cite{Ru10, RS06a, RS06b, Sh09, Sh12, Sh21, ShDR}.  One of the definitions of the nested CVaR is relevant to abstract dynamic programming \cite{Bert} dealing with equations $u_{N}(x)=\min_{a\in A(x)}H(x,a,u_{N-1}),$ for $x\in\X,$ where $u_N$ is an objective function for a finite or infinite time horizon $N,$  $x\in\X$ is the state,  $a\in A(x)$ is the action, and $u_0:\X\to\R$ is a terminal cost.   For the nested CVaR, the function $H$ is CVaR of $(c_a+\beta u_{N-1})$ with respect to a transition probability $p_a,$ where $c_a$ is the one-step cost.  For an MDP, the optimality equation for the nested CVaR for the risk level $\alpha \in [0,1]$ is
  \begin{equation}\label{eqURMDPCV}
u_N(x)=\min_{a\in A(x)}\max_{b\in B(x,\alpha,a)} \sum_{x'\in\X}[ c(x,a,x')+\beta u_{N-1}(x')]b_{x'}p(x'|x,a),\qquad x\in\X.
\end{equation}

Expressions \eqref{eqVRMDPCV} and \eqref{eqURMDPCV} are closely relevant.  In fact, \eqref{eqURMDPCV} can be viewed as \eqref{eqVRMDPCV} when the two-dimensional argument $(x,y)\in \X\times [0,1]$ of $v_N$ in \eqref{eqURMDPCV} is projected to a single-dimensional argument $x:=(x,\alpha).$  Also, similar to \eqref{eqURMDPCV}, formula \eqref{eqVRMDPCV} minimizes CVaR, but the tail risk level depends on the history of the process.  As algorithm DCVaR demonstrates, this risk level depends on the value function and previous gains and losses.  Contrary to this, for the nested CVaR, the risk level $\alpha$ is constant.  In this sense, the DCVaR is more flexible than the nested CVaR.

Let us define $V_N(x,y):=y v_N(x,y)$ for $N=0,1,\ldots$ or $N=\infty,$ $x\in\X,$ and $y\in [0,1]].$   Then \eqref{eqVRMDPCV} becomes
\begin{equation}\label{eqVRMDPCV_scaled}
V_N(x,y)=\min_{a\in A(x)}\max_{b\in B(x,y,a)} \sum_{x'\in\X}(yb_{x'}c(x,a,x')+\beta V_{N-1}(x',yb_{x'}))p(x'|x,a).
\end{equation}
This is the optimality equation for the DRMDP1 defined in the next paragraph.  It is easier to deal with the DRMDP1 than with the DRMDP because, as stated in Lemma~\ref{l71}, the value functions $V_N(x,y)$ are concave in $y$ when $V_0(x,y)$ is concave in $y\in [0,1],$ as assumed in Assumption~\ref{As1}.

 Let us introduce the DRMDP1 as an RMDP defined by the tuple $({\bf X},\A,\B,A(\cdot),B(\cdot,\cdot,\cdot),\tilde{c},q),$  where the state space $\bf{X},$ the action space $\A,$ uncertainty space $\B,$  sets of available actions $A(x),$ and uncertainty sets $B(x,y,a)$  are the same as in the DRMDP,  one-step costs
$\tilde{c}(x,y,a,b,x',y')=y'\cdot c(x,a,x'),$ and the transition probability
$q(x',D|x,y,a,b)=p(x'|x,a)\delta_{yb_{x'}}(D),$ where  $x,x'\in\X,$ $y\in [0,1],$ $D\in\mathcal{B}([0,1])$ $a\in A(x),$ and $b\in B(x,y,a).$   Let $V_N(x,y,\pi^\A,\pi^\B)$ be the expected total cost for the DRMD1.  
%
%
%
 The following theorem implies that the DRMDP and DRMDP1 correspond to the same problem when $\alpha>0,$ where $\alpha=y_0.$
 \begin{theorem}\label{teqMDMMDP1DD}
 For $x\in\X,$ $y\in [0,1],$ $\pi\in\Pi^\A,$   $\pi^\B\in\Pi^\B,$ and $N=0,1,\ldots,$ or $N=\infty,$
 \[
 V_N(x,y,\pi^\A,\pi^\B)=y v_N(x,y,\pi^\A,\pi^\B).
 \]
 \end{theorem}
 \begin{proof}
 The proof is based on the standard induction arguments. We write in the proof $x,y,a,b$ instead of $x_0,y_0,a_0,b_0,$ respectively. For $N=0$ the required equality holds since $V_0(x,y,\pi^\A,\pi^\B)=V_0(x,y)=yv_0(x,y).$  Let us show that it holds for horizon $N\ge 1$ if it holds for horizon $(N-1).$ Indeed, for $h_1^*=x,y,a,b$ and $y_1=yb_{x_1},$ by using \eqref{eqVRMDPCV} and notations from Appendix~\ref{a1},
\begin{eqnarray*}
 &\ & V_N(x,y,\pi^\A,\pi^\B)=\sum_{a\in A(x)}\pi_0^\A(a|x)\int_{b\in B(x,y,a)}\sum_{x_1\in \X}[yb_{x_1}c(x,a,x_1)\\
 &\ &+\beta yb_{x_1}v_{N-1}(x_1,yb_{x_1},\pi^{\A,h_1^*},\pi^{\A,h_1^*})p(x_1|x,a)\pi^\B_0(db|x,a)]=yv_N(x,y,\pi^\A,\pi^\B).
 \end{eqnarray*}
 Also, $V_\infty(x,u,\pi^\A,\pi^\B)=\lim_{N\to\infty}V_N(x,u,\pi^\A,\pi^\B)=\lim_{N\to\infty}yv_N(x,u,\pi^\A,\pi^\B)=yv_\infty(x,u,\pi^\A,\pi^\B).$
 \end{proof}
The following two corollaries follow directly from Theorem~\ref{teqMDMMDP1DD} and from formula  \eqref{eq4.2}, which holds in the form of equalities. 
\begin{corollary} \label{ceqMDMMDP1DD} For a finite or infinite horizon $N,$ initial state $x\in\X,$ and  initial risk level $y\in (0,1]$,  
the equality
${\rm DCVaR}_y(Z_N;x,\phi)=v_N(x,y) $ holds for a nonrandomized risk-independent policy $\phi\in\Pi$ for the DM if and only if  $V_N(x,y,\phi,\pi_*^\B)=V_N(x,y)$  for every optimal policy $\pi_*^\B\in\Pi_*^\B(N)$ for Nature.
\end{corollary}
\begin{corollary} \label{ceqMDMMDP1DD1} For a finite or infinite horizon $N,$ let $\phi\in \Pi$ be a nonrandomized risk-independent policy for the DM, $x\in\X,$ and $y\in (0,1].$ If   $V_N(x,y,\phi,\pi_*^\B)=V_N(x,y)$ for every optimal policy $\pi_*^\B \in \Pi_*^\B(N)$ for Nature, then
\begin{equation}\label{eq4.9}
{\rm DCVaR}_y(Z_N;x,\phi)=v_N(x,y)=\min_{\pi^\A\in\Pi^\A} {\rm DCVaR}_y(Z_N;x,\pi^\A)=\min_{\pi\in\Pi} {\rm DCVaR}_y(Z_N;x,\pi),
\end{equation}
and, in particular, $\phi$ is DM's optimal policy for the initial state $x.$
\end{corollary}


Formula \eqref{eqVRMDPCV_scaled} can be rewritten, for $N=1,2,\ldots$ and for $N=\infty,$ in the form of \eqref{eqfirstopb} and \eqref{eqsecondopa} as
\begin{equation}\label{eqVRMDPCVQ}
\begin{aligned}
Q_{N}(x,y,a)&=\max_{b\in B(x,y,a)} \sum_{x'\in\X}(yb_{x'}c(x,a,x')+\beta V_{N-1}(x',yb_{x'}))p(x'|x,a),\quad a\in A(x),\\
V_N(x,y)&=\min_{a\in A(x)}Q_{N}(x,y,a).
\end{aligned}\end{equation}
According to \eqref{eqAstar} and \eqref{eqBstar}, the sets of optimal actions for the DM are
\begin{equation}\label{eqAstarDCV}
A^*_{N}(x,y):=\bigl\{a\in A(x):\,V_{N}(x,y)=\max_{a\in A(x)} Q_N(x,y,a) \bigr\},\quad x\in{\bf X}, y\in[0,1],
\end{equation}
and for Nature, for  $x\in{\bf X},$ $y\in[0,1],$ $a\in A(x),$
\begin{equation*}\label{eqBstarDCV}\begin{aligned}
B^*_{N}(x,y,a):&=\bigl\{b^*\in B(x,y,a):\,Q_N(x,y,a)\\ &=\max_{b\in B(x,y,a)}  \sum_{x'\in\X}(yb_{x'}c(x,a,x')+\beta V_{N-1}(x',yb_{x'}))p(x'|x,a)\bigr\} .
\end{aligned}\end{equation*}

%
%
%


The following theorem provides a sufficient condition under which a policy for the MDP minimizes DCVaR.  The assumption that the policy $\phi\in\Pi\subset\Pi^\A$ is nonrandomized is not needed.  We formulate this assumption only because  this theorem is applied to a nonrandomized policies $\phi$ in this paper.
\begin{theorem}\label{tm4.5} For a nonrandomized risk-independent policy $\phi\in\Pi,$  initial state $x\in\X,$  tail risk level $\alpha\in (0,1],$ and time horizon $N=1,2,\ldots$ or $N=\infty,$ if $P_x^{\phi,\sigma^\B}(a_t\in A^*_{N-t}(x_t,y_t))=1,$
for every persistently optimal policy $\sigma^\B\in \Pi^\B$ for Nature and for every nonnegative integer $t<N,$ then formulae \eqref{eq4.9} hold.
\end{theorem}
\begin{proof}
 Let $\pi_*^\B\in\Pi^\B$ be an optimal policy for Nature. Then, in view of Corollary~\ref{corA2b}, there exists a persistently optimal policy $\sigma^\B\in\Pi^\B$ for Nature such that $V_N(x,\alpha,\phi,\pi_*^\B)=V_N(x,\alpha,\phi,\sigma^\B)=V_N(x,\alpha).$ Corollary \ref{ceqMDMMDP1DD1} implies \eqref{eq4.9}.
\end{proof}


If $\alpha=0,$ then \eqref{eqVRMDPCV} for $y=\alpha$ becomes
\begin{equation*}\label{eqVDET}
v_{N}(x,0)=\min_{a\in A(x)}\max_{x'\in\X} \{c(x,a,x')+\beta v_{N-1}(x',0): x'\in\X,\ p(x'|x,a)>0\}, \quad x\in\X,
\end{equation*}
which is the minimax equation for the sequential deterministic game in which the DM chooses actions, and Nature chooses a transition from the set of transitions having positive probability.  The DM is trying to minimize the length of a feasible path while Nature is trying to maximize it.   For a horizon $N=1,2,\ldots$ a path  is a sequence $x_0,a_0,x_1,a_1,\ldots,x_N,$ and for $N=\infty$ it is $x_0,a_0,x_1,a_1,\ldots.$ For a given initial state $x_0\in\X,$ a path is called feasible if $a_t\in A(x_t)$ and $p(x_{t+1}|x_t,a_t)>0$ for all integer values $t,$  $0\le t < N.$  The lengths of finite and infinite-horizon paths are defined by formulae \eqref{eqCN} and \eqref{eqC}  for finite and infinite horizon problems respectively. This is a special case of an RMDP considered in Appendix~\ref{a1}. In this case, the state space ${\bf X}:=\X$ is finite, action sets for the DM are $A(x),$ uncertainty sets for Nature are $B(x,a):=\{x'\in \X:\, p(x,a,x')>0\},$ one step costs are $c(x,a,x'),$ and all moves are deterministic to the next states $x'\in B(x,a)$  selected by Nature.  This problem can be solved easily by value iteration.

We notice that, if $\alpha=1,$ then all uncertainty sets can be reduced to singletons $B(x,a)=\{b^{x,a}\},$  where $b^{x,a}\in\R_+^M$ such that $b^{x,a}_{x'}=1$ if $p(x'|x,a)>0.$ Thus, if $\alpha=1,$ then we deal with an MDP with expected total discounted costs.

\section{Formulation of the Main Result}\label{s5}

Let us consider the DRMDP1.  We recall that $V_0(x,y)=yv_0(x,y),$ where the function $v_0$ is continuous in $y\in [0,1].$ When $N<\infty,$ everywhere in the rest of this paper, except Appendix~\ref{a1}, we assume that  Assumption~\ref{As1} holds. For example, this assumption holds when the function $v_0(x,y)$ is nonincreasing in $y$ and concave in $y$  for all $x\in\X.$ In particular, this assumption holds if the final payoff $v_0(x,y)$ does not depend on $y,$ that is, $v_0(x,y)=v_0(x).$
\begin{Assumption}\label{As1}
The function $V_0:\X\times [0,1]\to \R$ is concave in $y\in [0,1].$
\end{Assumption}

In this section we introduce Algorithm DCVaR constructing for $x\in\X$ and $\alpha\in (0,1]$ a nonrandomized policy $\phi\in\Pi$ minimizing ${\rm DCVaR}_\alpha(Z_N;x,\pi^\A)$ among all policies $\pi^\A\in\Pi^\A.$ In view of  \eqref{eqvalDVARDVAR}, this means   that  ${\rm DCVaR}_\alpha(Z_N;x,\phi)=v_N(x,\alpha).$      The case $\alpha=0$ is addressed in the previous section.  
Let $\alpha\in (0,1].$

%
%

We are interesting in dealing with concave functions because, according to Lemma~\ref{l71}, the functions $V_N(x,y)$ and $Q_N(x,y,a)$ are  concave in the tail risk level $y\in [0,1].$   In addition, if $N<\infty$ and $V_0(x,y)=yv_0(x),$ then the functions $V_N(x,y)$ and $Q_N(x,y,a)$ are piecewise linear in $y\in [0,1].$

 For  a real-valued continuous concave  function $y\mapsto f(y)$ defined on a finite  interval $[\bf{a},\bf{b}]\in\R,$ we denote by $f'^+(y)$ and   $f'^-(y)$ its right
 \[f'^+(y):=\lim_{\Delta y\downarrow 0} \frac{f(y+\Delta y)-f(y)}{\Delta y}
 \]
 and left
  \[f'^-(y):=\lim_{\Delta y\downarrow 0} \frac{f(y)-f(y-\Delta y)}{\Delta y}
 \]
 derivatives  respectively, where the concavity of $f$ implies $f'^+(y)\le f'^-(y).$ If $f$ is a real-valued  continuous concave function defined on the interval $[\bf{a},\bf{b}],$ we set  $f'^-({\bf a}):=+\infty$ and $f'^+({\bf b}):=-\infty.$ For a function of multiple variables, say $f(u,y,z),$ with all variables except one, say $y$, being discrete, we shall consider a function in $y$ for fixed values of all other parameters. The right and left derivatives in $y,$ if they exist,  will be also denoted as $f'^+(u,y,z) $ and $f'^- (u,y,z)$ respectively. If a derivative in $y$ exists for some $u,$ $y$, and $z,$ it is denoted by $ f'(u,y,z),$ which means that $f'(u,y,z):= f'^+(u,y,z)= f'^-(u,y,z).$
 We usually apply these notations to  functions of two variables $(x,y),$ where $x\in\X$ and $y\in [0,1].$  We recall that, for a real-valued continuous concave function $f$ defined on an interval,  the right and left derivatives  always exist, the function $f'^+$ is right-continuous and lower semicontinuous, and the function  $f'^-$ is left-continuous and upper semicontinuous. A concave function on an interval is differentiable everywhere except for at most a countable set. We also denote by $\partial_y f(\cdot,y,\cdot):=[f'^+(\cdot,y,\cdot), f'^-(\cdot,y,\cdot)]$   the superdifferential of $f$ in $y$ at the point $(\cdot,y,\cdot).$ In particular, if $f'(\cdot,y,\cdot)$ exists, then $\partial_y f(\cdot,y,\cdot):=\{f'(\cdot,y,\cdot)\}.$ We shall also apply this notation for the superdifferential of functions of one and two variables.


 For a real number $u\ge 0$ and a real-valued  continuous concave  function $f(y)$ defined on a finite interval $[\bf{a},\bf{b}],$  in view of the convention $f'^-(\bf{a}):=+\infty$ and $f'^ +(\bf{b}):=-\infty,$
 one of the following two possibilities holds: \begin{enumerate}[label=(\roman*)] \item   there is a unique point $y^*\in [\bf{a}, \bf{b}]$ such that $u\in\partial_y f(y^*);$
 \item there is an open interval $(\tilde{a},\tilde{b})\subset [\bf{a},\bf{b}]$ such that $f'(y)=u$ for all $y\in(\tilde{a},\tilde{b}).$
\end{enumerate}
In case (ii) there is a maximal open interval satisfying (ii).  If the function $f$ is piecewise linear, then $f'^-(y^*)> f'^+(y^*)$ in case (i).

Let us consider the sets of optimal actions $A^*_{t}(x,y)$ defined in \eqref{eqAstarDCV}. 
According to Theorem~\ref{tm4.5}, if for every optimal policy of Nature, at each finite epoch $t,$ such that $N>t\ge 0,$ a policy $\phi\in\Pi$ of the DM chooses  actions from the state $A^*_{N-t}(x_t,y_t),$ then this policy minimizes ${\rm DCVaR}_\alpha(Z_N;x,\pi^\A)$ in $\pi^\A\in\Pi^\A.$

The following equalities hold in view of Hiriart-Urruty and Lemar\'echal \cite[p.28]{HUL}:
\begin{equation}\label{partVQ+}
V'^+_t(x,y)=\min_{a\in A^*_t(x,y)} Q'^+_t(x,y,a),\quad x\in\X,\ y\in [0,1],\ t=1,2,\ldots,\infty,
\end{equation}
\begin{equation}\label{partVQ-}
V'^-_t(x,y)=\max_{a\in A^*_t(x,y)}  Q'^-_t(x,y,a),\quad x\in\X,\ y\in [0,1],\ t=1,2,\ldots,\infty.
\end{equation}
In particular, \eqref{partVQ+} and \eqref{partVQ-} imply that,
\begin{equation*}
\partial_y V_t(x,y)\subseteq \partial_y Q_t(x,y,a), \quad x\in\X,\ y\in [0,1],\ a\in A^*_t(x,y),\ t=1,2,\ldots,\infty.
\end{equation*}

{\bf Algorithm DCVaR} for running an optimal policy for an MDP with a finite state space $\X,$ finite action sets $A(z)$ available at states $z\in \X,$ transition probability $p,$ one-step costs $c,$  and discount factor $\beta\in (0,1],$ where $\beta<1 $ if $N=\infty.$

{\bf Inputs}: Initial state $x\in\X$ tail risk level $\alpha \in (0,1],$ the horizon length $N=1,2,\ldots$ or $N=\infty,$ the value functions $V_N(\cdot,\cdot), V_{N-1}(\cdot,\cdot),\ldots V_1(\cdot,\cdot),$ if $N<\infty,$ or the value function $V_\infty(\cdot,\cdot)$ if $N=\infty.$ After an action $a_t$ is selected at the state $x_t,$ where $t$ is a nonnegative integer such that $t<N,$ the next state $x_{t+1}$ such that $p(x_{t+1}|x_t,a_t)>0$ becomes known.

1. Set $x_0:=x,$ $y_0:=\alpha,$ $t:=0,$ $I:=1,$ and choose an arbitrary $\phi_t(x_t,y_t):=a_0\in A^*_{N}(x_0,y_0)$.

2. Do steps 2.1--2.4 while $t\le N-1$:

2.1. If $I=1,$ then choose an arbitrary number $u_{N-t}\in  [V_{N-t}'^+(x_t,y_t), V_{N-t}'^-(x_t,y_t)].$

2.2.  Compute:
\begin{equation}\label{eqkey}
u_{N-t-1}:=\frac{u_{N-t}-c(x_t,a_t,x_{t+1})}{\beta}.
\end{equation}

2.3. If there is a unique point $y^*\in [0,1]$ such that $u_{N-t-1}\in \partial_y V_{N-t-1}(x_{t+1},y^*),$ then set $y_{t+1}:=y^*$ and $I:=1.$
 Otherwise, set $I:=0$ and choose  $y_{t+1}\in (0,1)$ such that  $y_{t+1}$ is an interior point of the interval, on which the function $V_{N-t-1}(x_{t+1},\cdot)$ is linear with the slope $u_{N-t-1};$ in other words, choose     $y_{t+1}\in (0,1)$ such that                              $V'(x_{t+1},y_{t+1})=u_{N-t-1},$ and there are points
 $y^{(1)}\in (0,y_{t+1})$ and $y^{(2)}\in (y_{t+1},1)$ such that $ V_{N-t-1}'^+(x_{t+1},y^{(1)})=  V_{N-t-1}'^-(x_{t+1},y^{(2)})=u_{N-t-1}$.

2.4. If $t\ge N-1,$ then stop; otherwise,    choose an arbitrary optimal action $a_{t+1}\in A^*_{N-t-1}(x_{t+1},y_{t+1}),$ set   $\phi_{t+1}(x_{t+1},y_{t+1}),$ and then set $t:=t+1.$ 

\begin{figure}[!ht]
\centering
\begin{subfigure}{.5\textwidth}
  \centering
  \includegraphics[width=.95\textwidth]{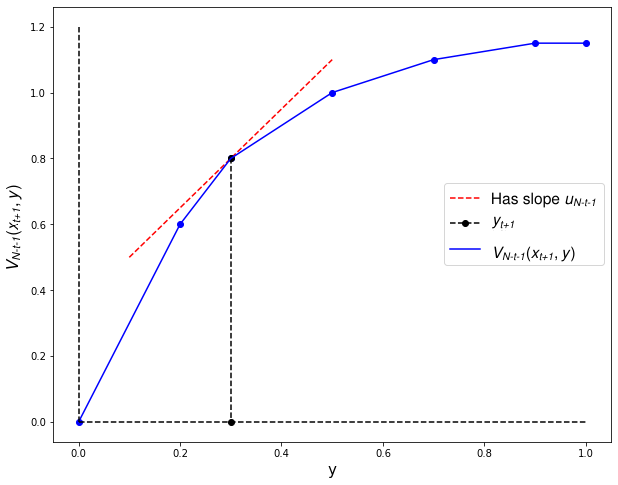}
  \caption{Case (i): $I=1$}
  \label{fig:demo1}
\end{subfigure}%
\begin{subfigure}{.5\textwidth}
  \centering
  \includegraphics[width=.95\textwidth]{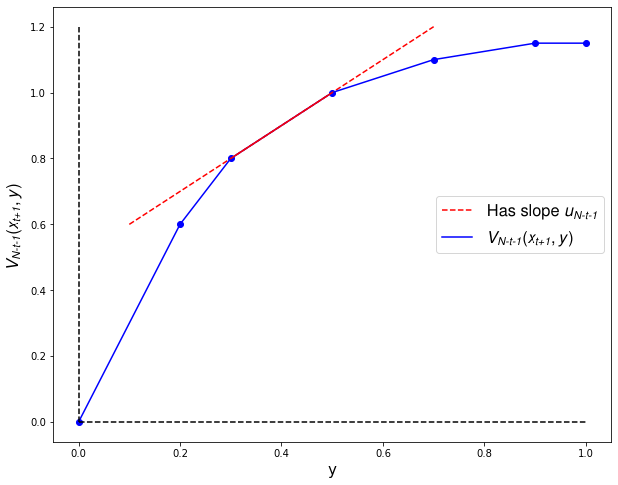}
  \caption{Case (ii): $I=0$}
  \label{fig:demo2}
\end{subfigure}
\caption{Demonstration of Cases $I=1$ and $I=0$ in Algorithm DCVaR for $t>0.$}
\label{fig:algo_demo}
\end{figure}

If $N<\infty$ then the algorithm stops after N iterations and returns a finite sentence of actions\\ $a_0,a_1,\ldots, a_{N-1}$ implementing an optimal policy at the steps $t=0,1,\ldots, N-1.$ If $N=\infty,$ then the algorithm returns an infinite sequence of actions $a_0,a_1,\ldots$ implementing an optimal policy at the steps $t=0,1,\ldots,$ and the algorithm can be stopped after a finite number of iterations because the impact of additional steps will be negligibly small due to the discount factor $\beta\in (0,1),$ and estimations of stopping times are standard.

 As an illustration of Algorithm DCVaR, Figure \ref{fig:algo_demo} shows the two main cases addressed by the algorithm. Figure \eqref{fig:demo1} is an example of a case where there is a unique tail risk level being identified in step 2.3, and Figure \eqref{fig:demo2} corresponds to the other case where there is an interval with the desired slope value.

The variable $I=1$ indicates that  the tail risk level $y_t$ is either given or identified  at step 2.3 as the unique point $y^*\in [0,1]$ such that $u_{N-t-1}\in\partial_y V_{N-t-1}(x_{t+1},y^*).$ The variable $I=0$ indicates that there is a nonempty open interval on which $u_{N-t-1}= V_{N-t-1}'(x_{t+1},y).$ If $N<\infty,$ then $y\mapsto V_{N-t-1}(x_{t+1},y)$ is a piecewise linear function and therefore the existence of the unique described point $y^*$ implies that $ V_{N-t-1}'^-(x_{t+1},y^*)> V_{N-t-1}'^+(x_{t+1},y^*).$
 If $N<\infty$ and $V_0(x,y)=yv_0(x),$ then the functions $V_n(x,\cdot)$ and $Q_n(x,\cdot,a)$  are piecewise linear for $n=1,2,\ldots,N,$ and  Subroutine 1 in Section~\ref{s6} computes functions $Q_n.$

Tail risk levels $y_t$ assigned by Nature are not available to the DM when $t>0.$  However, by using formula \eqref{eqkey}, which is based on the analysis of Nature's optimal policies,  at step 2.3 the algorithm  detects either the tail risk level $y_t=y^*$ or an  interval, in which $y_t$ is located, and every optimal action for an internal point of this interval is also optimal if Nature selects another tail risk level from this interval.
The main result of this paper is summarized in the following theorem, whose proof is provided in Section 8.

\begin{theorem}\label{TMAIN}
For $N=1,2,\ldots$ or $N=\infty,$ $x\in\X,$ and $\alpha\in (0,1],$ Algorithm DCVaR generates a nonrandomized risk-independent  policy $\phi\in\Pi$ minimizing DCVaR, that is,
\[
{\rm DCVaR}_\alpha(Z_N;P_x^\phi)=v_N(x,\alpha)=\min_{\pi^\A\in \Pi^\A}{\rm DCVaR}_\alpha(Z_N;P_x^{\pi^\A})=\min_{\pi\in \Pi}{\rm DCVaR}_\alpha(Z_N;P_x^\pi).
\]
\end{theorem}

\section{Properties of the Mass Transfer Problems Solved by Nature}\label{s6}

This section studies the problem Nature solves at each horizon $N=1,2,\ldots$ or $N=\infty.$  In particular, it studies problem \eqref{eqmcon} motivated by the maximization equation in \eqref{eqVRMDPCVQ}.  The results of this section are self-contained.  Continuity and concavity in $y\in [0,1]$ of the function $V(x,y)$ in formula \eqref{eqmcon} are assumed because, as shown in the next section, these assumptions are satisfied by the  functions $Q_N$ and $V_N$ from \eqref{eqVRMDPCVQ}.   Tneorem~\ref{p2} is the key fact on which formula \eqref{eqkey} in Algorithm  DCVaR is based.

A function $f:[a,b]\mapsto\R,$ where $a,b\in\R,$  is called piecewise linear, if there is a finite sequence of increasing numbers $(u_i)_{i=0}^n$ with $u_0=a$ and $u_n=b$ such that the function $f$ is affine on each interval $[u_{i-1},u_i],$ $i=1,\ldots,n.$

Let $V: \X\times[0,1] \mapsto \R^+$ be a real-valued function such that for each fixed $x\in\X=\{1,2,\ldots,M\} $ this function is  continuous and concave in $y\in [0,1].$  For a probability distribution $p$ on $\X,$ that is, $p(x)\ge 0$ for all $x\in\X$ and $\sum_{x\in\X} p(x)=1$, let
\begin{equation}\label{eqmcon}F(y): =\max_{b\in B(y)} \sum_{x\in\X} V(x,yb(x)) p(x),     \qquad\qquad             y\in [0,1],\end{equation}
where $B(y):=\{b\in\R^M: \  \sum_{x\in\X} p(x)b(x)=1,\ b(x)\ge 0,  yb(x)\le 1,\ x\in \X\}$.

A function $\tilde{b}:[0,1]\times\X \to\R^M$ is called feasible for problem \eqref{eqmcon} if for each $x\in \X$ the function $y\mapsto \tilde{b}(y,x)$ is Borel-measurable on $[0,1],$ and $\tilde{b}(y,\cdot)\in B(y)$ for all $y\in [0,1].$ For example, $\tilde{b}(y,x)\equiv 1$ is a feasible function. 
A feasible function $\tilde{b}(\cdot,\cdot)$ for problem \eqref{eqmcon} is called a  \emph{solution} to problem   \eqref{eqmcon} if $F(y)=\sum_{x\in\X} V(x,y\tilde{b}(y,x))p(x)$ for all $y\in [0,1].$

For each $y\in [0,1]$ we also consider the set $B^*(y)=\{b\in B(y): F(y)=\sum_{x\in\X} V(x,yb(x))p(x)\}$  of solutions at $y.$ In other words, a feasible function $\tilde{b}:[0,1]\times\X \to\R^M$ is a solution if and only if $\tilde{b}(y,\cdot)\in B^*(y)\}$ for all $y\in [0,1].$

The following theorem states the existence of solutions and  describes the properties of the  function $F.$
\begin{theorem}\label{p1}
 There exists a  solution to problem~\eqref{eqmcon}, and
  the function $F:[0,1]\to\R$  is  continuous   and concave.
\end{theorem}
 \begin{proof} Without loss of generality we assume that $p(x)>0$ for all  $x\in\X$ because, if  $p(x)=0,$ then state $x$ can be excluded from the set $\X.$
 We observe that $y\mapsto B(y)$ is a continuous compact-valued set-valued mapping of $[0,1]$ into the set of subsets of $\R^M.$  In addition $B(y)\subset B(0)$ for all $y\in [0,1],$ and the function $(y,b)\mapsto \sum_{x\in\X}V(x,yb(x))p(x)$ is continuous on $[0,1]\times B(0).$   Therefore, the Berge maximum theorem implies that  $B^*(y)$ are nonempty compact sets for all $y\in [0,1],$ the function $F$ is continuous, and the set-valued mapping $y\mapsto B^*(y)$ is upper-semicontinuous at all $y\in [0,1].$ In view of the Kuratowsli-Ryll-Nardzewski  measurable selection theorem, there is a Borel-measurable optimal solution $b,$ which implies that all functions $b(y,x)$ are Borel-measurable in $y\in [0,1].$ Let $y_1,y_2\in [0,1]$ and $y_3=\lambda y_1+ (1-\lambda)y_2,$ where $\lambda\in [0,1].$ Then $b^*(y_3)= \lambda b(y_1)+ (1-\lambda)y_2\in B(y_3).$  Then  $F(y_3)\ge  \sum_{x\in\X}V(x,y_3b^*_3(y_3,x))p(x)\ge \lambda F(y_1)+(1-\lambda)F(y_2),$ where the second inequality follows from concavity of $V.$
\end{proof}

We notice that concavity of $F:[0,1]\to\R$ implies  continuity of $F:(0,1)\to\R,$ but  this fact is not used in the proof of Theorem~\ref{p1}.
The following theorem describes the properties of the solutions to problem \eqref{eqmcon}.
\begin{theorem}\label{p2} Let $\tilde{b}$ be a  solution to problem~\eqref{eqmcon} and $x\in\X.$ If $p(x)>0,$ then for each $y\in [0,1]$
\begin{eqnarray}\label{ineqlevel}
 V'^-(x,y\tilde{b}(\tilde{y},x))&\ge& F'^-(y)\quad {\rm for}\ \tilde{y}\in [0,y],\label{eq6.02}\\
V'^+(x,y\tilde{b}(\tilde{y},x))&\le& F'^+(y)\quad {\rm for}\ \tilde{y}\in [y,1].\label{eq6.03}
\end{eqnarray}

In addition, the following relations hold for all $y\in [0,1]:$
\begin{eqnarray}\label{ineqkey}
\max \{ V'^+(x,y\tilde{b}(y,x)):\, p(x)>0,\,  x\in\X\}&\leq& \min \{V'^-(x,y\tilde{b}(y,x)):\, p(x)>0,\,  x\in\X  \},\label{eq6.2}\\
\label{eqkey1} F'^-(y)&=&\min\{ V'^-(x,y\tilde{b}(y,x)):p(x)>0, x\in\X\},\label{eq6.3}\\
\label{eqkey2} F'^+(y)&=&\max\{ V'^+(x,y\tilde{b}(y,x)):p(x)>0, x\in\X\},\label{eq6.4}
\end{eqnarray}
and, if  the function $V(x,y)$  is piecewise linear in $y$ for each $x\in \X, $ then the function $F(y)$ is also piecewise linear.
\end{theorem}
\begin{proof} Problem \eqref{eqmcon} can be simplified to a problem which has a natural interpretation.  First, without loss of generality we can assume that $V(x,0)=0.$   Indeed, we can always consider the objective function $\hat{V}(x,y):= V(x,y)-V(x,0).$ Then all the objective and value functions will be shifted by the constant $\sum_{x\in\X} V(x,0)p(x),$ and the solution sets won't change.
Second, in the rest of this section we assume  without loss of generality that $p(x)>0$ for all $x\in\X.$

Third, let us consider the functions $\tilde{z}(y,x):=y\tilde{b}(y,x)p(x),$ where $y\in [0,1],$ and vectors $z(x):=yb(x)p(x),$ where $b=(b(1),\ldots,b(M))\in\R^M,$  $x\in\X.$ For $x\in\X$ and $z\in [0,p(x)],$ let us define $v(x,z ):=
V(x,z/p(x))p(x).$ Then $v'^+(x,\tilde{z}(y,x))=V'^+(x,y\tilde{b}(y,x)),$  $v'^-(x,\tilde{z}(y,x))=V'^-x,y\tilde{b}(y,x)),$ and
problem  \eqref{eqmcon} becomes
\begin{equation}\label{eqmcon1}F(y): =\max_{z\in Z(y)} \sum_{x\in\X} v\left(x,z(x)\right),     \qquad\qquad             y\in [0,1],\end{equation}
where $Z(y):=\{z\in\R^M: \  \sum_{x\in\X} z(x)=y,\ 0\le z(x)\le p(x),\ x\in \X\}.$

Similarly to problem \eqref{eqmcon}, we can consider feasible functions and solutions $\tilde{z},$ and solution sets $Z^*(y):=\{z\in Z(y): F(y)=\sum_{x\in\X} v(x,z(y))p(x)\},$ $y\in [0,1],$ for problem \eqref{eqmcon1}.  Then $\tilde{b}(\cdot,\cdot)$ is a solution (or feasible function) for problem  \eqref{eqmcon} if and only if $\tilde{z}(y,x)=y\tilde{b}(y,x)p(x),$ where $y\in [0,1],$ $x\in\X,$ is a solution (a feasible function, respectively) for problem \eqref{eqmcon1}.

Thus, in order to prove Theorem~\ref{p1}, it is sufficient to prove its statement for problem~\eqref{eqmcon1} with the functions $V$ and $\tilde{b}$ replaced with $v$ and $\tilde{z}$ respectively.  In particular,  formulae \eqref{ineqkey}--\eqref{eqkey2} become
\begin{eqnarray}\label{eqkeya}
\max \{ v'^+(x,\tilde{z}(y,x)):\,   x\in\X\}&\leq& \min \{v'^-(x,\tilde{z}(y,x)):\,   x\in\X  \},\\
\label{eqkey1a} F'^-(y)&=&\min\{v'^-(x,\tilde{z}(y,x)):\, x\in\X\},\\
\label{eqkey2a} F'^+(y)&=&\max\{v'^+(x,\tilde{z}(y,x)):\, x\in\X\}.
\end{eqnarray}

Problem \eqref{eqmcon1} describes the following optimal mass transfer problem.  For  $x=1,\ldots,M,$ let us consider  $M$ intervals $W(x):=\{x\}\times [0,p(x)],$ $x\in\X,$  in $\R^2,$ which we call sources. These sources can be interpreted as vertical cylindrical vessels of the same diameter filled with a liquid up to  heights $p(x).$ Let $W:=\cup_{x\in \X} W(x).$ There is an additional
 interval $ [0,1],$ which we call the destination. This interval can be interpreted as an empty vertical cylindrical vessel with the same diameter  and with  height 1. The total value of the liquid at source $x$ from level 0 up to the level $z\in [0,1]$ is $v(x,z),$ where $v(x,\cdot)$ is a concave function with $v(x,0)=0,$ $x\in\X.$
 The goal is to fill up the destination by the liquid from sources to maximize  for each $y\in [0,1]$ the total value of the liquid in the destination from the level 0 up to the level $y$.  Any amount of the liquid can be taken from any part of each source.

We recall, that for a concave function $f:[a,b]\to\R$ with  $-\infty<a<b<+\infty$ and $f(a)=0,$ for $y\in[a,b]$
 \[
 f(y) = \int_a^y f'^+(z)dz=\int_a^y f'^-(z)dz,
 \]
  $f'^+(y)\le f'^-(y),$ and the set $\{y\in [a,b]: f'_+(y)<f'_-(y)\}$ is countable. Thus, both functions $v'^+(x,y)$ and $v'^(x,y),$ where $y\in [0,p(x)],$ describe the nonlinear unit cost of the liquid at the source $x\in\X$ depending on the height  $y.$ These unit costs may  have arbitrary signs (positive, negative, or 0), but they are nondecreasing functions of the hight $y\in [0,p(x)]$ since the function $v(x,y)$ is concave in $y.$

  Let us prove \eqref{eq6.02} and \eqref{eq6.3}. If $y=0,$ then these formulae hold in the form $+\infty=+\infty.$   Let $y\in (0,1].$  Since the function $v(x,\cdot)$ is   concave, the most valuable part of the liquid of the volume $y\in [0,p(x)]$ at the source $x\in\X$ is the interval $[0,y],$ whose value is $v(x,y).$ Thus, if the amounts $z(y,x)$ of liquid should be taken by an optimal policy for the level $y$ from sources $ x\in\X,$ then the decision to take all the liquid from the interval $[0,z(y,x)]$ from each source $x\in\X$ is optimal, and $v'^-(x,\tilde{y})\ge F'^-(y)$  for all $\tilde{y}\in [0,z(y,x)].$ 
  In particular, $v'^-(x,z(\tilde{y},x))\ge F'^-(\tilde{y})\ge F'^-(y)$ for $\tilde{y}\in [0,y].$  So. \eqref{eq6.02} is proved.

Inequality  \eqref{eq6.02} can be written as  $\min_{x\in\X}v'^-(x,z(\tilde{y},x))\ge F'^-(y)$ for $\tilde{y}\in (0,y]$ if $y\in(0,1].$ In particular, $\min_{x\in\X}v'^-(x,z(y,x))\ge F'^-(y),$ and the strict inequality is impossible.  Indeed, if the strict inequality holds, then concavity of $v(x,\cdot)$ and $F(\cdot)$ imply the  existence of $\epsilon\in (0,y)$ such that $\min_{x\in\X}v'^-(x,z(y_1,x))> F'^-(y_1)$ for all $y_1\in [y-\epsilon,y].$
  Since $\sum_{x\in\X}z(y-\epsilon,x)=y-\epsilon<y= \sum_{x\in\X}z(y,x),$ there exists $x^*\in\X$ such that $z(y-\epsilon,x^*)<z(y,x^*).$ Let $d:=z(y,x^*)-z(y-\epsilon,x^*).$ Then we consider the feasible solution with the same decisions outside of the interval of $[y-\epsilon,y-\epsilon+d]$ and allocating all possible resource from the source $x^*$ when the level is in this interval.  The formal definition is: for $y_2\in [0,1]$
     \[z_1(y_2,x)=\begin{cases}z(y_2,x), &\text{if $y_2\notin [y-\epsilon,y-\epsilon +d]$ and $x\in\X;$}\\
                                                                     z(y-\epsilon,x), &\text{if $y_2\in [y-\epsilon,y-\epsilon +d]$ and $x\in\X\setminus \{x^*\}$;}\\
                                                                     z(y-\epsilon,x^*)+ y_2-(y-\epsilon),& \text{if $y_2\in [y-\epsilon,y-\epsilon +d]$ and $x= x^*.$}
  \end{cases}\]
  Let us verify that $z_1$ is a feasible function for problem \eqref{eqmcon1}.  We need to verify the properties of $z_1(y_2,x)$ for $y_2\in [y-\epsilon,y-\epsilon +d].$  Let $y_2\in [y-\epsilon,y-\epsilon +d].$ Then $0\le z_1(y_2,x^*)=  z(y-\epsilon,x^*)+ y_2-y+\epsilon\le  z(y-\epsilon,x^*)+d=   z(y,x^*)\le p(x^*).$  In addition,
  $\sum_{x\in\X}z_1(y_2, x)=y_2.$ Thus, $z_1$ is feasible, and
  \begin{equation}\begin{array}{ll}\label{eq6fin}
 & F(y-\epsilon+d)= F(y-\epsilon)+\int_{y-\epsilon}^{y-\epsilon+d}F'^-(y_2)dy_2\\ &< F(y-\epsilon)+\int_{y-\epsilon}^{y-\epsilon+d}v'^-(x^*,z_1(y_2,x^*))dy_2 =\sum_{x\in\X}v(x,z_1(y-\epsilon+d,x)),
\end{array}\end{equation}
which is impossible since $F$ is an optimal value, where the last equality in \eqref{eq6fin} follows from
\[
F(y-\epsilon)=\int_0^{y-\epsilon}\sum_{x\in\X}v'^-(x,z_1(y_2,x))dy_2.
\]
Formula \eqref{eqkey1a} is proved, and it implies \eqref{eqkey1}.


  The proof of \eqref{eq6.4} is similar. In particular, $b(1,x)=1$ for all $x\in\X$ is the only solution for $y=1,$ and \eqref{eq6.4} holds in the form of $-\infty=-\infty.$  So, we consider $y\in [0,1).$  The same arguments as for \eqref{eq6.02} imply \eqref{eq6.03}.  Therefore, $F'^+(y)\ge \max\{v'^+(x,\tilde{z}(y,x)):\, x\in\X\},$ where $z$ is a solution of \eqref{eqmcon1} corresponding to a solution $b$ of \eqref{eqmcon}. The similar argument as in the proof of \eqref{eq6.3} imply that the strict inequality is impossible because otherwise there is a solution $\tilde{z}$ such that $\sum_{x\in\X}v(\tilde{z}(\tilde{y},x)>F(\tilde{y})$ for some $y\in (y,1).$ Thus, \eqref{eqkey2a} and therefore \eqref{eq6.4} hold. Equations \eqref{eq6.3} and \eqref{eq6.4} imply \eqref{eq6.2}. The last claim is correct because, if functions $V(x,y)$ are piecewise linear in $y$ for each $x\in\X$, then formulae \eqref{eq6.3} and \eqref{eq6.4} imply that the  functions $F'^-(y)$ and $F'^+(y)$ take finite numbers of values when $y\in [0,1].$ Therefore, the concave function $F$ is piecewise linear on the interval $[0,1].$
\end{proof}
\begin{remark}{\rm
It is possible to construct all solutions of problem \eqref{eqmcon}.  We do not use these solutions in this paper.  So, we do not provide details on constructing them. If $y=0$ then $z(x)\equiv0$ is the unique solution of problem \eqref{eqmcon1}. If $y=1,$ then $z(x)=p(x),$ $x\in\X,$ is the unique solution  of \eqref{eqmcon1}. Thus, $b(x)\equiv 1$ is a solution of problem \eqref{eqmcon} for $y\in\{0,1\},$ and for $y=1$ this solution is unique. Let $y\in (0,1).$ Then, let us consider two cases. In case 1, $y$ does not belongs to an  interval $(y_*,y^*)\subset [0,1]$ on which the concave function $F$ is linear.     In case 2, $y$  belongs to an  interval $(y_*,y^*)\subset [0,1]$ on which the function $F$ is linear.  In case 2 we choose $(y_*,y^*)$ in the way that this is the maximal open interval containing $y$ on which $F$ is linear. In case 1, there is a unique solution $b(y,x),$ $x\in\X,$ at $y$. In case 2,  $X(y)\ne\emptyset,$ where $X(y)$ is the set of all states $x\in\X$ such that there are points $y_*(x),y^*(x)\in [0,1]$ such that $y_*(x)<y^*(x),$ and the function $v(x,\cdot)$ is linear on $(y_*(x),y^*(x))$ with its derivative on this interval equal to  $F' (y).$   In this case we consider $(y_*(x),y^*(x))$ being the maximal open interval satisfying this condition.    If $y,y_1\in (y_*,y^*)$ then $X(y)=X(y_1).$   In case 2, $y^*-y_*=\sum_{x\in\X(y)}(y^*(x)-y_*(x))p(x),$ and $b(x,\cdot)$ is a solution at $y$ iff: (i) $b(y,x)=b(y_*(x),x)$ for $x\notin X(y),$ and (ii) $y-y_*=\sum_{x\in\X(y)}(b(y,x)-y_*(x))p(x)$ and $b(y,x)\in [y_*(x),y^*(x)]$ for $x\in X(y).$ We notice that, in case 2, if $X(y)=\{x^*\}$ is a singleton, then  $b(y,x):=b(y_*,x),$ if $x\in\X\setminus \{x^*\},$ and $b(y,x^*):=b(y_*,x^*)+(y-y_*),$ is the unique solution of problem \eqref{eqmcon} at $y\in (0,1).$}
\end{remark}
\begin{remark} {\rm
In order for solutions to have physical meaning of moving the liquid from sources to the destination, the optimal functions  $z(y,x)$ should be nondecreasing in $y.$ To achieve this, it is sufficient  to define appropriately the values  $z(y,x)$ for $y\in (0,1)$ when $y$ belongs to an interval $(y_*,y^*)$ on which the function $F(y)$ is linear, and $X(y)$ is not a singleton.   In this case, for example, we can move the liquid from sources $x\in X(y)$ to the destination  sequentially.  Sources with smaller numbers can be used first.  Such solutions are nondecreasing in $y.$ Of course, there are other ways to construct nondecreasing solutions.}
\end{remark}

 A concave piecewise linear function $f: [a,b]\to\R $ on a finite interval can be represented by a finite sequence $\{(q^f_i,l_i^f)\}_{i=1,\ldots,I^f},$ where $q_i^f$ are slopes (derivatives at linear intervals) and $l_i^f$ are lengths of linear intervals $i=1,\ldots,I^f,$ where   $q^f_i>q^f_{i+1},$ $l^f_i>0,$ and $\sum_{i=1}^{I^f} l^f_i=b-a.$ Formally speaking, $f'(y)=q_i^f,$ if $a+\sum_{j=1}^{i-1}l^f_j <y < a+\sum_{j=1}^i l^f_j ,$ where $i=1,\ldots,I^f$ and $\sum_{i=1}^0:=0.$

For $x\in\X,$ if a function $V(x,y)$ is represented by a sequence $\{(q^{V(x,\cdot)}_i,l_i^{V(x,\cdot)})\}_{i=1,\ldots,I^{V(x,\cdot)}},$ then the function $v(x,z)$ is represented by the sequence $\{(q^{v(x,\cdot)}_i,l_i^{v(x,\cdot)})\}_{i=1,\ldots,I^{v(x,\cdot)}}$ with $q^{v(x,\cdot)}_i=q^{V(x,\cdot)}_i,$ $l^{v(x,\cdot)}_i=p(x)l^{V(x,\cdot)}_i ,$ and $I^{v(x,\cdot)}=I^{V(x,\cdot)}.$ In view of Theorem~\ref{p2}, if all the functions $V(x,y)$ are piecewise linear in $y,$ then the function $F$ is piecewise linear, and it can be constructed in the following way.

%

\noindent\textbf{Subroutine 1:}

1. Merge sequences $\{(q^{v(x,\cdot)}_i,l_i^{v(x,\cdot)})\}_{i=1,\ldots,I^{v(x,\cdot)}},$ $x\in\X,$ into a single finite sequence  $\{(q_i,l_i)\}.$
By doing this, merge the intervals with the same slopes.  This means that any final set of pairs $(q,l_{i_j})$ with distinct indexes $i_j$ and the same slope $l_{j_i},$  where $j=1,\ldots,J$ and $J\le M,$ should be replaced with the single  pair $(q,\sum_{j=1}^J l_{i_j}).$

3. The resulted finite sequence $\{(q^F_i,l_i^F)\}_{i=1,\ldots,I^F},$ where $\max_{x\in\X} I^{v(x,\cdot)} \le I^F\le\sum_{x\in\X} I^{v(x,\cdot)}$ and $q^F_i>q^F_{i+1},$ represents the concave piecewise linear  function $F.$

Subroutine 1 is the major step in recursive computations of value functions $V_N(\cdot,\cdot), Q_N(\cdot,\cdot,
\cdot)$ for $N<\infty;$ see the following section.

\section{Properties of Value Functions and Sets of Optimal Actions.}\label{s7}

This section describes the properties of the value functions $Q_N(x,y,a)$ and $V_N(x,y)$ for DRMDP1 defined in equations \eqref{eqVRMDPCVQ}, where $N=1,2,\ldots$ or $N=\infty.$ It also describes the properties of the sets of optimal actions for the DM.  We recall that the function $V_0(x,y)=yv_0(x,y)$ is concave and continuous in $y$ because the function $v_0(x,y)$ is continuous in $y,$ and because Assumption~\ref{As1} holds. We observe that, if $v_0(x,y)=v_0(x),$ then the function $V_0(x,y)=yv_0(x)$ is linear in $y\in[0,1].$

We observe that for fixed $x\in\X,$ $a\in A(x),$ and $N=1,2,\ldots$ or $N=\infty,$ the first equation in \eqref{eqVRMDPCVQ} is a particular instance of equation \eqref{eqmcon} with $F(y):= Q_N(x,a,y),$ $B(y):=B(x,y,a),$ with the variable $x$ in \eqref{eqmcon} substituted with $x',$  with $V(x',y):=yc(c,a,x')+\beta V_N(x',y), $ and with $p(x'):=p(x'|x,a).$

\begin{lemma}\label{l71}
For $N=1,2,\ldots,$ and for $N=\infty,$ the functions $Q_N(x,y,a)$ and $V_{N}(x,y)$ defined in  \eqref{eqVRMDPCVQ}  are continuous and concave in $y\in[0,1]$ for all $x\in\X$ and $a\in A(x).$ Furthermore, if $N<\infty$ and $V_0(x,y)$ is piecewise linear in $y\in [0,1],$ then these functions are piecewise linear in $y\in [0,1].$
\end{lemma}
\begin{proof}
For $N=1,2,\ldots$ we prove this lemma by induction.
 For $N=0$, the function $V_0(x,y)$ is continuous and concave in $y.$ This is explained in the first paragraph of this section.

Assume  that the function $V_N(x,y)$ is continuous and concave in $y\in [0,1]$ for some $N=0,1,\ldots.$
Theorem~\ref{p1}  applied to the first formula in \eqref{eqVRMDPCVQ} implies that the functions $Q_{N+1}(x,y,a)$ is continuous and concave in $y,$ and, in view of the second formula in \eqref{eqVRMDPCVQ},      the function $V_{N+1}(x,y)$ is continuous and concave in $y.$ If the function $V_0(x,y)$ is piecewise linear in $y\in [0,1],$ then the last claim in Theorem~\ref{p2} applied to the first formula in \eqref{eqVRMDPCVQ} implies that the functions $Q_{N+1}(x,y,a)$ is piecewise linear in $y.$ Since for each $x\in\X$  the set $A(x)$ is finite, the second formula in   \eqref{eqVRMDPCVQ} implies that the functions $V_{N+1}(x,y,a)$ are piecewise linear in $y.$

In view of \eqref{equnbanch},  concave and continuous in $y$ functions $V_N(x,y)$ converge  uniformly to $V_\infty(x,y).$ Thus, $V_\infty(x,y)$ is concave and continuous in $y.$  
\end{proof}

Thus, for each $x\in\X$ and for each $N=0,1,\ldots$ or $N=\infty,$ the function $V_N(x,y)$ is concave and continuous in $y$ on the interval $[0,1].$  We recall that by definition $V'^-_N(x,0)=+\infty$ and $V'^+_N(x,1)=-\infty$. Then for each real number $d$ exactly one of the following to possibilities takes place: (i) either there exists a unique $y\in [0,1]$ such that $d\in \partial_{y}V_N(x,y)$, or (ii) there exist a unique interval $[{\bf a},{\bf b}]\subset [0,1]$ such that on this interval the function $V_N(x,\cdot)$ is linear with the slope $d$, and $d\notin \partial_{y}V_N(x,y)$ if $y\in [0,{\bf a})\cup ({\bf{b}},1].$

\begin{lemma}\label{l72}
Assume that, for $N=1,2,\ldots$ or $N=\infty$ and for some $x\in\X,$ the value function $V_N(x,y)$ is linear on an interval  $[{\bf a},{\bf b}],$ where $0\le {\bf a}< {\bf b}\le 1.$ Then the following statements hold:  \begin{enumerate}[label=(\roman*)] \item
 if  $a\in A^*_N(x,y)$ for some $y\in ({\bf a},{\bf b}),$ then $a\in A^*_N(x,y)$ for all $y\in [{\bf a},{\bf b}];$
\item  $A^*_N(x,y)=A^*_N(x,\tilde{y})$ for  $y,\tilde{y}\in ({\bf a},{\bf b});$
\item $ Q'_N(x,y,a)= V'_N(x,y)$ for $y\in ({\bf a},{\bf b})$ and $a\in A^*_N(x,y).$
\end{enumerate}
\end{lemma}

\begin{proof}
According to Lemma~\ref{l71}  the functions $V_N(x,y)$ and $Q_N(x,y,a),$ $a\in A(x),$ are concave in $y\in [0,1],$ and
 $V_N(x,y)\le Q_N(x,y,a)$ for all $y\in [0,1].$
 Consider statement (i) and let $a\in A^*_N(x,\tilde{y})$ for some $\tilde{y}\in ({\bf a},{\bf b}).$  This means that  $V_N(x,\tilde{y})= Q_N(x,\tilde{y},a).$ Since $V_N$ and $Q_N$ are concave functions, and $Q_N$ dominates $V_N,$ we have that   $V_N(x,y)= Q_N(x,y,a)$ for all $y\in  [{\bf a},{\bf b}].$ Statements (ii) and (iii) follow from (i).
\end{proof}

The state space ${\bf X}:= \X\times [0,1]$ is infinite, and functions $Q_N(x,y,a)$ and $V_N(x,y)$ can be computed by discretizing the interval $[0,1].$  They can be also computed exactly if the state and action sets $X,\  A$ and the horizon length $N$ are not too large.  For each $N=1,2,\ldots,$ computing $Q_N(x,y,a)$ can be done by Subroutine 1 described at the end of Section~\ref{s6}. According to the second equation in \eqref{eqVRMDPCVQ}, for each $x\in\X$ functions  $V_N(x,y)$ are lower envelops in $ a\in A(x)$ of piecewise-linear concave in $y\in [0,1]$ functions $ Q_N(x,y,a),$  and their computation is a standard task in computational geometry.

\section{Proof of Theorem~\ref{TMAIN}.}\label{s8}
This section contains the proof of Theorem~\ref{TMAIN}.

\begin{proof}[Proof of Theorem~\ref{TMAIN}] For the given time horizon $N=1,2,\ldots$ or $N=\infty$ and for a given sequence $x_0,x_1,\ldots,$ Algorithm DCVaR sequentially generates a finite or infinite sequence of optimal actions.  Let $x_0,x_1,$ $x_2,\ldots$ be the states of the system and  $y_0=\alpha,y_1,y_2,\ldots$ be the tail risk levels, and the values $y_t$ are not observed by the DM when $t\ge 1.$  In addition to the parameters of the model, the algorithm also uses the value functions $V_N(\cdot,\cdot), V_{N-1}(\cdot,\cdot),\ldots, V_{N-1}(\cdot,\cdot)$, if $N<\infty,$ or it uses the value function $V_\infty(\cdot,\cdot)$ if $ N=\infty.$

The parameter $I$ takes two values: 1 and 0.  $I=1$ means that at the current time epoch $t$ the tail risk level $y_t$ is known, and $I=0$ indicates that it is not known. However, in the second case the algorithm finds an interval of possible risk levels, and this is sufficient to choose an optimal action.

At each iteration $t=0,1,\ldots,$ the algorithm chooses an action $a_t$. In view of Theorem~\ref{tm4.5}, in order to prove the theorem, it is sufficient to show that $a_t\in A^*_{N-t}(x_t,y_t)$ for all nonnegative integer $t< N$ for every persistently optimal policy played by Nature.

The proof is based on induction in $t.$ According to step 1,  $a_0\in A^*_N(x_0,y_0).$ Let an optimal action $a_t\in A^*_{N-t}(x_t,y_t)$ be selected at an epoch $t=0,1,\ldots$ by step 1 or 2.4 of the algorithm. In view of \eqref{eqVRMDPCVQ},
\[
V_{N-t}(x_t,y_t)=Q_{N-t}(x_t,y_t,a_t)=\max_{b\in B(x,y,a)} \sum_{x'\in\X}(yb_{x'}c(x,a,x')+\beta V_{N-1}(x',yb_{x'}))p(x'|x,a).
\]
  Formula \eqref{eqkey2} implies
\begin{equation*}\begin{array}{ll}
Q'^+_{N-t}(x_t,y_t,a_t)&=\max\{c(x_t,a_t,x')+\beta V'^+_{N-t-1}(x',y_tb_{x'}):\,x'\in\X,\ p(x'|x_t,a_t)>0\}\\
&\ge  c(x_t,a_t,x_{t+1})+\beta V'^+_{N-t-1}(x_{t+1},y_{t+1}),
\end{array}\end{equation*}
where the last inequality holds because $x_{t+1}\in\X,$ and $y_{t+1}=y_tb_{x_{t+1}}$ for every $b\in B^*(x_t,y_t,a_t).$ Thus,
\begin{equation}\label{eq8.1}
V'^+_{N-t-1}(x_{t+1},y_{t+1})\le\frac{ Q'^+_{N-t}(x_t,y_t,a_t)-c(x_t,a_t,x_{t+1})}{\beta}.
\end{equation}
Similarly to \eqref{eq8.1}, formula \eqref{eqkey1} implies
\begin{equation}\label{eq8.2}
\frac{ Q'^-_{N-t}(x_t,y_t,a_t)-c(x_t,a_t,x_{t+1})}{\beta}\le    V'^-_{N-t-1}(x_{t+1},y_{t+1}).
\end{equation}

Therefore, in view of inequalities \eqref{eq8.1} and \eqref{eq8.2},
\begin{equation*}\label{eqkey88}
u_{N-t-1}:=\frac{u_{N-t}-c(x_t,a_t,x_{t+1})}{\beta}\in [V'^+_{N-t-1}(x_{t+1},y_{t+1}), V'^-_{N-t-1}(x_{t+1},y_{t+1})]=\partial_y V_{N-t-1}(x_{t+1},y_{t+1})
\end{equation*}
for $u_{N-t}\in \partial_y V_{N-t}(x_t,y_t).$  If there is a unique point $y^*$ such that $u_{N-t-1}\in \partial_y V_{N-t-1}(x_{t+1},y^*),$ then  $y_{t+1}=y^*$ is the tail risk level at the state $x_{t+1}$ and epoch $(t+1).$

If there are multiple points $y^*$ such that $u_{N-t-1}\in \partial_y V_{N-t-1}(x_{t+1},y^*),$ that is, the case $I=0$ takes place, then concavity of $V_{N-t-1}(x_{t+1},\cdot)$  implies that there is a maximal interval $[{\bf a},{\bf b}]\subset [0,1]$ such that the function  $V_{N-t-1}(x_{t+1},y)$ is linear in $\in [{\bf a},{\bf b}],$ and its slope is $u_{N-t-1}.$ According to Lemma~\ref{l72}, for all $y\in ({\bf a},{\bf b})$ the sets $A^*_{N-t-1}(x_{t+1},y)$ coincide,  and $ Q'_{N-t-1}(x_{t+1},y,a)=V'_N(x_{t+1},y)$ for  $a\in A^*_{N-t-1}(x_{t+1},y).$ Therefore, step 2.1 can be skipped, and future calculations do not depend on the particular value $y_{t+1}\in ({\bf a},{\bf b}).$  So, $a_t\in A^*_{N-t}(x_t,y_t)$ in spite of the fact that the DM does not know the states $y_t$ when $t>0.$  The algorithm either calculates values $y_t$ or can calculate intervals to which values $y_t$ belong, and optimal action sets coincides for  tail risk levels, which are internal points these intervals. Actions from these set are also optimal for the extreme points ${\bf a}$ and ${\bf b}$ of this interval. \end{proof}


\section{Extension to Stochastic Cost Functions}\label{s9}
 In this section, we extend the results of this paper to random one-step cost functions with finite supports.  Arbitrary random one-step costs can be approximated by such costs.  

  In practice, one-step costs $c(x,a,x')$ can be random. To model random costs, for $x,x'\in\X$ and $a\in A(x),$ let us consider finite sets $W(x,a,x')$ and  real-valued  one-step costs $c(x,a,x',w'),$ where $w'\in W(x,a,x').$ If a control $a$ is selected at a state $x,$ and the system moves to a state $x',$ a random cost $c(x,a,x',w')$ is collected,
  where $w'$ has a discrete distribution $q(\cdot|x,a,x')$ satisfying $q(w'|x,a,x')\geq 0$ for all $w'\in W(x,a,x')$, and $\sum_{w'\in W(x,a,x')} q(w'|x,a,x') = 1.$

  Let us set $\W:=\cup_{x,x'\in\X,a\in A(x)} W(x,a,x').$ Let $c(x,a,x',w'):=0$ and $q(w'|x,a,x'):=0$ for $w'\in \W\setminus W(x,a,x').$

  If we augment the state  $x$ with the parameter $w\in \W,$ we have a particular case of the original problem with the expanded state space.
 To explain details, let the augmented state be $(x,w),$ where $w\in\W$ is the realization of the random outcome that took place at the previous time instance.  At the initial time 0, there are no previous events.  In this case,  we choose an arbitrary $w_0\in \W$ and consider an initial state $(x_0,w_0)$ instead of the original initial state $x_0.$

  So, we obtain the MDP with the state space $\X\times\W,$ action space $\A,$ sets of available actions $A(x,w):=A(x),$ transition probabilities
 \[\tilde{p}(x',w'|x,w,a) = p(x'|x,a)q(w'|x,a,x'),\quad x,x'\in\X,a\in A(x),w,w'\in \W,\]
 and one-step costs
 \[\tilde{c}(x,w,a,x',w') = c(x,a,x',w'),\quad x,x'\in\X,a\in A(x),w,w'\in \W.\]
Notice that the transition probabilities and costs do not depend on $w.$ Therefore, the value functions also do not depend on the state component $w.$  In view of this detail, we provide below the minimax equations for DRMDP1 for this model, which are similar to \eqref{eqVRMDPCV_scaled}:

\begin{equation}\label{ereqQ_w}
\begin{aligned}
\tilde{Q}_{N+1}(x,y,a) = \max_{b\in \tilde{B}(x,y,a)}\sum_{(x',w')\in\X\times\W}  \left(y b_{(x',w')} c(x,a,x',w')+\beta\tilde{V}_N(x',y b_{(x',w')})\right)q(w'|x,a,x')p(x'|x,a),
\end{aligned}
\end{equation}
and
\begin{equation}\label{ereqoa_w}
\tilde{V}_{N+1}(x,y)=\min_{a\in A(x)}\tilde{Q}_{N+1}(x,y,a),
\end{equation}
 where for all $x\in\X,$ $a\in A(x),$ and $y\in[0,1]$
\begin{equation*}\label{BCVaR_w}
\begin{aligned}
\tilde{B}(x,y,a)&=\{b\in R^{\tilde{M}}:\, 0\leq yb_{(x',w')}\leq 1,\forall x'\in\X,w'\in W(x,a,x'), \\
&\sum_{(x',w')\in\X\times\W} b_{(x',w')}p(x'|x,a)q(w'|x,a,x') = 1,\ {\rm and}\ b_{(x',w')}=0\ {\rm if}\ p(x'|x,a)q(w'|x,a,x')=0\\
&{\rm or}\ w'\notin W(x,a,x') \}
\end{aligned}
\end{equation*}
with $\tilde{M}: = M|\W|,$ where $M$ and $|\W|$ are the numbers of elements of the sets $\X$ and $\W$ respectively.

As a remark, minimax equations \eqref{ereqQ_w} and \eqref{ereqoa_w} can be extended to possibly infinite sets $\X,$ $\A,$ and $\W.$  For metric spaces $\X,$ $\A,$ $\W$ and nonnegative costs $c,$ equation
 \eqref{ereqQ_w} can be rewritten in the integral form
\begin{equation*}\label{ereqQ_wcont}
\begin{aligned}
&\tilde{Q}_{N+1}(x,y,a) = \max_{b\in \bar{B}(x,y,a)}\{\int_\X\int_\W (y b(x',w') c(x,a,x',w')+\beta\tilde{V}_n(x',y b(x',w')))q(dw'|x,a,x')p(dx'|x,a)\},
\end{aligned}
\end{equation*}
where
\begin{equation}\label{BCVaR_wcont}
\begin{aligned}
\bar{B}(x,y,a)&=\{b\in\mathbb{F}:\, 0\leq yb(x',w')\leq 1,\forall x'\in\X,w'\in W(x,a,x'),\\ &\int_\X\int_\W b(x',w')q(dw'|x,a,x')p(dx'|x,a) = 1,\ {\rm and}\ b(x',w')=0\ {\rm if}\ w'\notin W(x',a,x)\},
\end{aligned}
\end{equation}
and $\mathbb{F}$ is a set of measurable real-valued functions on $\X\times\W.$  Developing specific conditions for correctness of \eqref{BCVaR_wcont} for infinite sets $\X,$ $\Y,$ and $\W$ is beyond the scope of this paper.

\section*{Acknowledgement}
This research was partially supported by the U.S. Office of Naval Research (ONR) under Grant\\ N000142412608.
The authors thank Jefferson Huang and Pavlo Kasyanov for their comments and express special thanks to Alexander
Shapiro for valuable suggestions and insightful discussions.

\begin{appendices}
\section{Results on Robust MDPs Used in this Paper}\label{a1}

 This appendix provides results on RMDPs used in this paper. In general, for robust optimization problems, some parameters of the model are unknown, and the goal is to select the best solutions for the worst possible values of unknown parameters.   Starting from \cite{GI05} and \cite{NG} there are numerous studies of RMDPs motivated by applications.  For an RMDP, there is an additional parameter on which one-step costs and transition probabilities depend, and at each step Nature (sometimes called Player II) selects this parameter from a set of available parameters.  The DM (sometines called Player I) chooses actions by knowing the past and current states and the past actions chosen by the DM and by Nature. Nature chooses actions by knowing this information and, in addition, by knowing the current action chosen by the DM.   Such models were also studied under the names of two-person zero-sum games with perfect information (see e.g., \cite{AFS, Gi, JN}) and minimax control \cite{GHH02}. In RMDP papers \cite{GI05, NG} the unknown parameter is the transition probability, and one-step costs depend only on states and actions selected by the DM.  However, this specification,  probably motivated by applications, does not simplify the  theory needed in this paper.  The models studied in \cite{GI05, NG} are currently called the models satisfying the property of $(s,a)$-rectangularity. Another property called $s$-rectangularity was introduced in \cite{WKR}, and there are other rectangularity properties; see \cite{GD, LiSh, GGC, WSBZ}. Models with $s$-rectangularity are relevant to two-person zero-sum stochastic games \cite{JN, Shap}, when the second player does not know the current action selected by the first player.

 For the purposes of this paper, it is sufficient to consider an RMDP with a state space being a compact subset of a Euclidean space, with a finite action space for the DM,   with action sets of Nature being compact subsets of a Euclidean space, with continuous one-step cost functions, and with weakly continuous transition probabilities.    We  consider  a slightly more general model in this appendix.

%

 An RMDP is defined by a tuple $({\bf X},\A,\B,A(\cdot),B(\cdot,\cdot),c,q).$ Here ${\bf X}$ is the state space, and there are two players: the DM and Nature.  There are two decision sets $\A$ and $\B,$ where $\A$ is the set of actions for the DM, and $\B$ is the set of actions for Nature.  We assume that ${\bf X},$ $\A$, and $\B$ are nonempty Borel subsets of Polish (complete, separable, metric) spaces. 
 The sets $A(x)\subset\A,$ where $x\in{\bf X},$ are assumed to be nonempty and finite, and $A(x)$ are the sets  of actions available to the DM at states $x\in{\bf X}.$
 For each pair $(x,a),$ where $x\in{\bf X}$ and $a\in A(x),$ a set of feasible actions $B(x,a)$ of Nature is defined.  We assume that the sets $B(x,a)$ are nonempty and compact subsets of $\B$ for all $x\in{\bf X}$ and for all $a\in A(x).$   We also assume that  $\Gr\ A:=\{(x,a):\, x\in{\bf X}, a\in A(x)\}$ is a Borel subset of  ${\bf X}\times\A,$ and the set-valued mapping  $B:\Gr\ A \mapsto 2^\B$ is continuous. Continuity of $B$ implies Borel measurability of $\Gr\ B:=\{(x,a,b):\, (x,a)\in \Gr\ A, b\in B(x,a)\}.$ According to the Arsenin-Kunugui measurable selection theorem, each of the graphs $\Gr\ A$ and $\Gr\ B$ contains at least one Borel measurable function. This fact, which is needed for the existence of policies, also follows from the Kuratowski and Ryll-Nardzewski measurable selection theorem.

The one-step cost $c(x,a,b,x')$ depends on the current state $x\in{\bf X},$ action $a\in A(x)$ chosen by the DM, action $b\in B(x,a)$ chosen by Nature, and the next state $x'\in {\bf X}.$ The function $c: ({\rm Gr}\ B)\times{\bf X}\to\R$ is assumed to be bounded and continuous.   The transition probability $q$ from ${\rm Gr}\ B$ to ${\bf X}$ is weakly continuous in $(x,a,b)\in{\rm Gr}\ B.$ Thus, $x'\sim q(\cdot|x,a,b),$ and $\int_{\bf X} f(x')q(dx'|x,a,b)$ is a continuous function on ${\rm Gr}\ B$ for every bounded continuous function $f:{\bf X}\to\R.$ 

 At each time step $t=0,1,\ldots$ Nature knows the action chosen by the DM  at time $t$ when Nature makes a decision.  Let $H_t:=({\bf X}\times\A\times\B)^t\times{\bf X}$ be the sets of histories that can be observed by the DM at time  $t.$ The set of histories, that can be observed by Nature at time $t,$ is $H_t\times \A.$ A policy $\pi^\A$ of the DM is a sequence
 $(\pi^\A_t)_{t=0,1,\ldots}$ of regular transition probabilities from $H_t$ to $\A$ such that $\pi^\A_t(A(x_t)|x_0,a_0,b_0,\ldots,x_{t-1},a_{t-1},b_{t-1},x_t)=1.$ Nature's policy $\pi^\B$   is a sequence
 $(\pi^\B_t)_{t=0,1,\ldots}$ of regular transition probabilities from $H_t\times\A$ to $\B$ such that $\pi^\B_t(B(x_t,a_t)|x_0,a_0,b_0,\ldots,x_{t-1},a_{t-1},b_{t-1},x_t,a_t)=1.$ Let $\Pi^\A$ and $\Pi^\B$ be the sets of policies for the DM and  Nature respectively.    Similarly to MDPs, an initial state $x\in{\bf X}$ and a pair of policies $(\pi^\A,\pi^\B)\in\Pi^\A\times\Pi^\B$ define the probability  $P_x^{\pi^\A,\pi^\B}$ on the space of trajectories $H:=({\bf X}\times\A\times\B)^\infty.$  An expectation with respect to this probability is denoted by $\E_x^{\pi^\A,\pi^\B}.$

 It is also possible to consider nonrandomized, Markov, and deterministic policies. A policy $\pi^\A$ for the DM is called nonrandomized if for each $t=0,1,\ldots$ there exists a mapping $\phi^\A_t:H_t\to\A$ such that  $\pi^\A_t(\phi^\A_t(h_t)|h_t)=1$ for all $h_t=x_0,a_0,b_0,\ldots,x_{t-1},a_{t-1},b_{t-1},x_t\in H_t.$ Such a policy is also denoted by $\phi^\A.$  A Markov policy $\phi^\A$ for the DM is a sequence of measurable functions $\phi_t^\A:{\bf X}\to\A,$ $t=0,1,\ldots,$ such that $\phi_t^\A(z)\in A(z)$ for all $z\in{\bf X}.$   A deterministic policy for the DM  is a Markov policy $\phi^\A$ such that $\phi^\A_t(z)=\phi^\A_s(z)$ for all $t,s=0,1,\ldots$ and all $z\in\bf{X}.$

 A policy $\pi^\B$ for Nature is called nonrandomized, if for each $t=0,1,\ldots$ there exists a mapping $\phi^\B_t:H_t\times\A\to\B$ such that  $\pi^\B_t(\phi^\B_t(h_t,a_t)|h_t,a_t)=1$ for all $h_t=x_0,a_0,b_0,\ldots,x_{t-1},a_{t-1},b_{t-1},x_t,a_t\in H_t\times\A.$ Such a policy is also denoted by $\phi^\B.$  A Markov policy $\phi^\B$ for Nature is a sequence of measurable functions $\phi_t^\A:{\bf X}\times\A\to\B,$ $t=0,1,\ldots,$ such that $\phi_t^\B(z,a)\in B(z,a)$ for all $(z,a)\in{\bf X}\times\A.$   A deterministic policy for Nature  is a Markov policy $\phi^\B$ such that $\phi^\B_t(z,a)=\phi^\B_s(z,a)$ for all $t,s=0,1,\ldots$ and all $(z,a)\in\bf{X}\times\A.$

If an initial state is $x\in{\bf X}$ and the DM and Nature play policies $\pi^\A$ and $\pi^\B$ respectively, then the expected total finite-horizon payoff of the DM to Nature is
\[
  v_N(x,\pi^\A,\pi^\B):=\E_x^{\pi^\A,\pi^\B}[ \sum_{t=0}^{N-1}\beta^t c(x_t,a_t,b_t,x_{t+1})+v_0(x_N)],\quad N=1,2,\ldots,\ \beta\in [0,1],
\]
where $v_0:{\bf X}\mapsto\R$ is a bounded continuous function (the terminal payoff), and for the infinite horizon
\[
  v_\infty(x,\pi^\A,\pi^\B):=\E_x^{\pi^\A,\pi^\B}\left[ \sum_{t=0}^{\infty}\beta^t c(x,a_t,b_t,x_{t+1})\right],\quad \beta\in [0,1).
\]

 For  $\beta\in [0,1]$ and $N=1,2,\ldots$  let us define sequentially
 \begin{equation}\label{eqfirstopb} Q^*_N(x,a):=\max_{b\in B(x,a)}\int_{\bf{X}} [c(x,a,b,x')+ \beta v_{N-1}(x')]q(dx'|x,a,b),\quad x\in {\bf X},\ a\in A(x),\end{equation}
 \begin{equation}\label{eqsecondopa}
 v_N(x):=\min_{a\in A(x)} Q^*_N(x,a),\quad x\in {\bf X},
 \end{equation}
where, in view of Berge's maximum theorem, each function $Q^*_N(x,a)$ is continuous,
 and therefore $v_N(x)$ is a continuous function.
 Equations \eqref{eqfirstopb} and \eqref{eqsecondopa} imply the minimax equation
\begin{equation}\label{eqminimaxgejn}
v_N(x)=\min_{a\in\A(x)}\max_{b\in B(x,a)} \int_{\bf X}[c(x,a,b,x')+ \beta v_{N-1}(x')]q(dx'|x,a,b), \quad x\in{\bf X}.
\end{equation}

Let $\beta\in [0,1).$  Then all continuous functions $v_N$ are uniformly bounded by $C/(1-\beta),$ where $C= \max\{\sup\{|v_0(x)|:\,x\in {\bf X}\},\sup\{|c(x,a,b,x')|:\, (x,a,b)\in {\rm Gr}\ B,\ x'\in{\bf X}\}\}.$  Banach's fixed point theorem applied to the space of uniformly bounded continuous functions on $\bf X$ implies that the minimax operator applied in \eqref{eqminimaxgejn} to the function $v_{N-1}$  has a unique fixed point $v_\infty(x),$ which is a bounded continuous function,   and $v_N(x)\to v_\infty(x)$ as $N\to\infty.$  In particular
\begin{equation}\label{equnbanch}
\sup_{x\in\bf{X}} |v_N(x)-v_\infty(x)|\le C\beta^N/(1-\beta)\to 0\ {\rm as}\ N\to\infty.
\end{equation}
Let $Q^*_\infty(x,a)$ be defined in \eqref{eqfirstopb} with $N=\infty.$ Thus,  \eqref{eqsecondopa} holds also for $N=\infty,$ and
$v_\infty$ is a unique bounded continuous function satisfying \eqref{eqminimaxgejn} or \eqref{eqfirstopb}, \eqref{eqsecondopa} with $N=\infty.$

For  $N=1,2,\ldots$ and for $N=\infty,$ let us consider the following nonempty sets of actions for the DM
\begin{equation}\label{eqAstar}
A^*_{N}(x):=\bigl\{a\in A(x):\,v_{N}(x)= Q^*_N(x,a) \bigr\},
\end{equation}
minimizing the right-hand side of \eqref{eqsecondopa}, and for Nature, for $x\in{\bf X}$ and  $a\in A(x),$
\begin{equation}\label{eqBstar}
B^*_{N}(x,a):=\bigl\{b^*\in B(x,a):\,Q^*_N(x,a) = \int_{\bf X}[\{c(x,a,b^*,x')+ \beta v_{N-1}(x')]q(dx'|x,a,b^*) .
\end{equation} maximizing the right-hand side of \eqref{eqfirstopb}.
The sets $A_N^*(x)$ are finite. As follows from the Berge maximum theorem, the sets $B_N^*(x,a)$ are compact.

For an $N$-horizon problem with $N=1,2,\dots$ or $N=\infty,$ a pair of policies $(\pi^\A_*,\pi^\B_*)\in \Pi^\A\times\Pi^\B$ is called  an \emph{equilibrium} if for all policies $\pi^\A\in\Pi^\A$ and $\pi^\B\in\Pi^\B$ 
\begin{equation*}\label{eqA7}
v_N(x,\pi^\A_*,\pi^\B)\le v_N(x,\pi^\A_*,\pi^\B_*)\le v_N(x,\pi^\A,\pi^\B_*),\qquad x\in{\bf X}.
\end{equation*}
\begin{theorem}\label{thmA0} For a given horizon $N=1,2,\ldots$ or $N=\infty,$ if two pairs of policies $(\pi_1^\A,\pi_1^\B)$  and $(\pi_2^\A,\pi_2^\B)$  are equilibria, then $(\pi_1^\A,\pi_2^\B)$  and $(\pi_2^\A,\pi_1^\B)$ are equilibria too, and $v_N(x,\pi_i^\A,\pi_j^\B)=v_N(x,\pi_1^\A,\pi_1^\B)$ for all $x\in {\bf X}$ and $i,j=1,2.$\end{theorem}
\begin{proof} The definition of an equilibrium implies that $v_N(x,\pi_1^\A,\pi_1^\B)=v_N(x,\pi_2^\A,\pi_2^\B)$ for  all $x\in {\bf X}.$ This equality and the definition of equilibria imply the conclusions of the theorem.
\end{proof}
If a pair $(\pi^\A_*,\pi^\B_*)\in \Pi^\A\times\Pi^\B$ is an equilibrium for some $\pi^\B_*,$ then $\pi^\A_*$ is called an optimal policy for the DM, and   $\pi^\B_*$ is called an optimal policy  for Nature.


In order to formulate structural properties of optimal policies,  we need to introduce additional notations and definitions. Let us fix a horizon $N=1,2,\ldots$ or $N=\infty.$


Let the DM and Nature reset the clock to 0 at time $t=1,2,\ldots.$  Let the history prior to time $t$ be $h^*_t=x^*_0,a^*_0,b^*_0,\ldots,x^*_{t-1},a^*_{t-1},b^*_{t-1}\in H_{t-1}\times \A\times\B$, and let $h_s=x_0,a_0,b_0,\ldots,x_{s-1},a_{s-1},b_{s-1},x_s\in H_s$ and $h_s,a_s$ be the histories the DM and Nature observe respectively during the following  $s=0,1,\ldots$  units of time after the clock is reset to 0.   If the clock is not  not reset, then at time $(t+s)$ the history would be $\tilde{h}_{t+s}=h^*_t,h_s$ for the DM and   $\tilde{h}_{t+s},a_s=h^*_t,h_s,a_s$ for Nature.

Let $\pi^{\A,h^*_t}$ and  $\pi^{\B,h^*_t}$ the policies for the DM and for Nature, if the clock is reset to 0 at time $t=1,2,\ldots$  and the finite sequence $h^*_t$  of past states and actions is observed. For $s=0,1,\ldots$
\[\pi_s^{\A,h_t^*}(a_s|h_s):=\pi_{t+s}^\A (a_s|h^*_t,h_s)\quad {\rm and} \quad \pi_s^{\B,h_t^*}(db_s|h_s,a_s):=\pi_{t+s}^\A (db_s|h^*_t,h_s,a_s)\]
are the policies the DM and Nature will use respectively following their original policies $\pi^\A$ and $\pi^\B,$ if the clock is reset to 0 at time $t.$

\begin{definition}A policy $\pi^\A$ for the DM ($\pi^\B$ for Nature) is called \emph{perfectly optimal} for an $N$-horizon problem, if it is optimal and, for each positive integer $t<N$ and for each  $h^*_t=x^*_0,a^*_0,b^*_0,\ldots,x^*_{t-1},a^*_{t-1},b^*_{t-1}\in H_{t-1}\times \A\times\B,$ the policy $\pi^{\A,h^*_t}$ ($\pi^{\B,h^*_t}$) is optimal for the horizon $(N-t).$
\end{definition}
\begin{definition}
A policy $\pi^\B$ for Nature is called \emph{persistently optimal} for an $N$-horizon problem, if\\ $\pi_t^\B( B^*_{N-t}(x_t,a)|h_t,a)=1$ for all nonnegative integers $t<N,$  all $h_t=x_0,a_0,b_0,\ldots, x_{t-1},a_{t-1},b_{t-1},x_t\in H_t,$ and all $a\in A(x_t).$
\end{definition}
We have that a nonrandomized policy $\phi^\B$ for Nature is persistently optimal if and only if $\phi_t(h_t,a)\in B^*_{N-t}(h_t,a)$ for all nonnegative integers $t<N,$  all $h_t=x_0,a_0,b_0,\ldots, x_{t-1},a_{t-1},b_{t-1},x_t\in H_t,$ and all $a\in A(x_t).$ A Markov  policy $\phi^\B$ for Nature is persistently optimal if and only if $\phi_t^\B(x,a)\in B^*_{N-t}(x,a)$ for all nonnegative integers $t<N,$  all $x\in {\bf X},$ and all $a\in A(x).$  For $N=\infty$ a deterministic optimal policy $\phi^\B$ for Nature is persistently optimal if and only if $\phi(x,a)\in B^*_\infty(x,a)$ for   all $x\in{\bf X},$ and all $a\in A(x).$
\begin{theorem} \label{thmA1}
For every horizon $N=1,2,\ldots$ or $N=\infty,$ the following statements hold:

   $(i_A)$ a policy $\pi^\A$ for the DM  is perfectly optimal  if and only if
$\pi_t^\A( A^*_{N-t}(x_t)|h_t)=1$  for all nonnegative integers $t<N$ and for all $h_t=x_0,a_0,b_0,\ldots, x_{t-1},a_{t-1},b_{t-1},x_t\in H_t;$

$(i_B)$ a persistently optimal policy $\pi^\B$ for Nature is perfectly optimal;

$(ii_A)$ there exists a Markov perfectly optimal policy for the DM, and a Markov policy $\phi^\A$ for the DM   is perfectly optimal if and only if $\phi^\A_t(x)\in A^*_{N-t}(x)$  for all nonnegative integers $t<N-1$ and for all $x\in {\bf X};$

$(ii_B)$ there exists a Markov persistently optimal policy for Nature, and  $\phi^\B$ is a Markov persistently optimal policy for Nature if and only if $\phi^\B_t(x,a)\in B^*_{N-t}(x,a)$ for all nonnegative integers $t<N-1,$ all $x\in {\bf X},$ and all $a\in A(x);$

$(iii)$ if $\pi^\A$ and $\pi^\B$ are optimal policies for the DM and Nature respectively, then
\begin{equation}\label{eqvNxxAB}
v_N(x)=v_N(x,\pi^\A,\pi^\B),\qquad x\in {\bf X};
\end{equation}

$(iv_A)$ for $N=\infty$ there exist deterministic optimal policies for the DM and Nature, and all deterministic optimal policy $\phi^\A$ for the DM and $\phi^\B$ for Nature   are perfectly optimal;

$(iv_B)$ for $N=\infty,$ a deterministic policy $\phi^\A$ for the DM is optimal if and only if $\phi(x)\in A_\infty(x)$  for all $x\in {\bf X},$ and every deterministic persistently optimal policy $\phi^\B$ for Nature is optimal.
\end{theorem}

\begin{proof}  Let $N<\infty.$  Let $\pi^\A$ be a policy for the DM satisfying the conditions stated in $(i_A),$  and let $\pi^\B$ be a persistently optimal policy for the Nature. Let us prove that the policies $\pi^\A$ and $\pi^\B$ are perfectly optimal for the DM and Nature respectively.  For a history $h_t=x_0,a_0,b_0,\ldots, x_{t-1},a_{t-1},b_{t-1},x_t,$ let us denote by $h_t^*$ its part that does not include the last state, that is, $h^*_t:=x_0,a_0,b_0,\ldots, x_{t-1},a_{t-1},b_{t-1}.$ We also define $v_0(x_N,\sigma^{\A,h^*},\sigma^{\B,h_N^*}):=v_0(x_N)$ for all policies. Straightforward backward calculations imply that if $v_t(x_{N-t},\pi^{\A,h_{N-t}^*},\pi^{\B,h_{N-t}^*})=v_t(x_{N-t})$ for $t=0,1,\ldots,N-1$ and for all $h_{N-t}\in H_{N-t},$ then  $v_{t+1}(x_{N-t-1},\pi^{\A,h_{N-t-1}^*},\pi^{\A,h_{N-t-1}^*})=v_{t+1}(x_{N-t-1})$  for all $h_{N-t-1}\in H_{N-t-1}.$  In particular, \eqref{eqvNxxAB} holds when $t=N-1.$
If $\sigma^\A\in\Pi^\A$ and $\sigma^\B\in\Pi^\B,$ then backward calculations imply that
\[
v_t(x_{N-t},\pi^{\A,h_{N-t}^*},\sigma^{\B,h_{N-t}^*}) \le   v_t(x_{N-t}) \le   v_t(x_{N-t},\sigma^{\A,h_{N-t}^*},\pi^{\B,h_{N-t}^*})\  {\rm if}\ t=1,\ldots,N, h_{N-t}\in H_{N-t}.
\]
Thus, $\pi^\A$ and $\pi^\B$ are perfectly optimal. So, statement ${\rm (i_B)}$ and the sufficient condition in ${\rm (i_A)}$ are proved. Since \eqref{eqvNxxAB} holds for persistently optimal policies $\pi^\A$ and $\pi^\B$ considered above,  Theorem~\ref{thmA0} implies that \eqref{eqvNxxAB} holds for any pair of equilibrium policies.  Thus (iii) is proved.

Let us prove the necessary condition in ${\rm (i_A)}.$  Let $\pi^\A$ be a persistently optimal policy for the DM, and it does not satisfy the conditions stated in ${\rm (i_A)}.$ Let $\pi^\B$ be a persistently optimal policy for Nature.  In view of ${\rm (i_B)}$, $\pi^\B$ be a perfectly optimal policy for Nature.  Let us consider the largest nonnegative integer $t< N$ such that $\pi_t^\A( A^*_{N-t}(x_t)|h_t)<1$ for some $h_t\in H_t.$ Then backward calculations imply that $ v_t(x_{N-t},\pi^{\A,h_{N-t}^*},\pi^{\B,h_{N-t}^*})> v_t(x_t),$ and thus the policy $\pi^\A$ is not perfectly optimal. This contradiction implies the necessary condition in ${\rm (i_A)}.$  Statement ${\rm (ii_A)}$ follows from ${\rm (i_A)}$ when $\pi^\A$ is a Markov policy.  Statement ${\rm (ii_B)}$ follows from definitions of persistently optimal policies and Markov policies for Nature.

Thus, statements ${\rm (i_A, i_B, ii_A, ii_B, iii)}$ are proved for $N<\infty.$ Let us prove them for $N=\infty.$ To prove ${\rm (i_A, i_B)},$ let us consider a policy $\pi^\A$ for the DM such that $\pi_t^\A(A_\infty^*(x_t)|h_t)=1,$ for all $t=0,1,\ldots$ and all $h_t\in H_t.$ We also choose a persistently optimal policy $\pi^\B$ for Nature for $N=\infty.$   Let us consider a finite-horizon problem with the terminal value $v_0=v_\infty.$ Then finite-horizon values for this problem with horizon $n=1,2,\ldots$ are equal to $v_\infty,$ and the sets $A_n^*(x)$ and $B_n^*(x,a)$ for this problem are equal to $A_\infty^*(x)$ and $B_\infty^*(x,a)$ respectively.  Statements  ${\rm (i_A})$ for a finite horizon implies that $\pi^\A$ is perfectly optimal for the DM, and $\pi^\B$ is perfectly optimal for Nature for all $n=1,2,\ldots.$  However, $\beta^n v_\infty (x)\to 0$ uniformly in $x\in{\bf X}$, and this implies that policies $\pi^\A$ and $\pi^\B$ are perfectly optimal for the DM and Nature respectively for $N=\infty.$ In particular, statement (iii) holds. To prove the necessary condition in ${\rm (i_A}),$ let $\pi^\A$ be a perfectly optimal policy for the DM.  Then $\pi^\A$ is perfectly optimal for every finite-horizon problem described in this paragraph.  Thus, $\pi_t^\A(A_\infty^*(x_t)|h_t)=1,$ for all $t=0,1,\ldots$ and all $h_t\in H_t.$  In view of Theorem~\ref{thmA0}, statements ${\rm (i_A})$ and ${\rm (ii_A})$ imply ${\rm (ii_A, ii_B, iii, iv_A, iv_B)}.$ In particular, a deterministic policy $\phi^\A$ for the DM )is perfectly optimal if and only if $\phi(x)\in A_\infty^*(x)$ for all $x\in{\bf X}.$  A deterministic policy $\phi^\B$ for Nature is persistently optimal if and only if $\phi(x,a)\in B_\infty^*(x)$ for all $x\in{\bf X},$ $a\in A(x).$
\end{proof}

%
%
\begin{corollary}\label{corA2b} For $N=1,2,\ldots$ or $N=\infty,$ if $\pi^\B\in\Pi^\B$ is an optimal policy for Nature, $x\in\X,$ and $\pi^\A\in\Pi^\A$ is a policy for the DM satisfying the condition
	$P_x^{\pi^\A,\pi^\B} (a_t\in A^*_{N-t}(x_t))=1$ for all nonnegative integers $t<N,$ then the following statements hold:
	
	(i) $v_N(x, \pi^\A,\pi^\B)=v_N(x);$
	
	(ii) $P_x^{\pi^\A,\pi^\B} (b_t\in B^*_{N-t}(x_t,a_t))=1$ for all nonnegative integers $t<N,$
	and there exists a persistently optimal policy for Nature $\sigma^\B\in\Pi^\B$ such that  $P_x^{\pi^\A,\sigma^\B} (a_t\in A^*_{N-t}(x_t))=1$ for all nonnegative integers $t<N,$ and   $v_N(x, \pi^\A,\sigma^\B)=v_N(x).$
\end{corollary}
\begin{proof} (i) Let $\tilde{\pi}^\A$ be a perfectly optimal policy for the DM. Then, in view of Theorem~\ref{thmA1}, $\tilde{\pi}^\A_t(A^*_{N-t}(x_t)|h_t)=1$  for all $h_t=x_0,a_0,b_0,\ldots, x_{t-1},a_{t-1},b_{t-1},x_t\in H_t,$ $t<N,$ and $v_N(x,\tilde{\pi}^\A,\pi^\B)=v_N(x).$

Let us define the perfectly optimal policy $\sigma^\A$  for the DM: for nonnegative integer $t<N$ and for $h_t\in H_t,$
\[
\sigma^\A_t(\cdot|h_t)=
\begin{cases}
\pi^\A_t(\cdot|h_t), &\text{if $\pi^\A_t(A^*_{N-t}(x_t)|h_t)=1,$}\\
\tilde{\pi}^\A_t(\cdot|h_t) &\text{otherwise}.
\end{cases}
\]
Since $P_x^{\pi^\A,\pi^\B} (a_t\in A^*_{N-t}(x_t))=1$  for all nonnegative integer $t<N,$ the probability measures   $P_x^{\pi^\A,\pi^\B}$ and $P_x^{\sigma^\A,\pi^\B}$  coincide. Therefore, $v_N(x,\pi^\A,\pi^\B)=v_N(x,\sigma^\A,\pi^\B)=v_N(x).$

(ii) Let $\tilde{\pi}^\B$ be a persistently optimal policy for Nature. Then, in view of Theorem~\ref{thmA1}, $\tilde{\pi}^\B_t(B^*_{N-t}(x_t,a_t)|h_t,a_t)=1$  for all $h_t=x_0,a_0,b_0,\ldots, x_{t-1},a_{t-1},b_{t-1},x_t\in H_t,$ $a_t\in A(x_t),$ $t<N,$

Let us define the persistently optimal policy $\sigma^\B$  for Nature: for nonnegative integer $t<N,$  for $h_t\in H_t,$ and for $a_t\in A(x_t),$
\[
\sigma^\B_t(\cdot|h_t,a_t)=
\begin{cases}
\pi^\B_t(\cdot|h_t,a_t), &\text{if $\pi^\B_t(B^*_{N-t}(x_t,a_t)|h_t,a_t)=1,$}\\
\tilde{\pi}^\B_t(\cdot|h_t,a_t) &\text{otherwise}.
\end{cases}
\]

Let assume that it is possible that $P_x^{\pi^\A,\pi^\B} (b_t\in B^*_{N-t}(x_t,a_t))<1$ for some nonnegative integer $t<0.$ Then
\begin{equation*}\label{eqbad}
P_x^{\sigma^\A,\pi^\B} (b_t\in B^*_{N-t}(x_t,a_t))=P_x^{\pi^\A,\pi^\B} (b_t\in B^*_{N-t}(x_t,a_t))<1.
\end{equation*}
Thus, while the DM plays the perfectly optimal policy $\sigma^\A,$ the optimal policy $\pi^\B$ for Nature can be improved on the set of that has a positive probability if Nature switches to the policy $\sigma^\B.$ Standard arguments imply that   $v_N(x,\sigma^\A,\sigma^\B)>v_N(x,\sigma^\A,\pi^\B)=v_N(x),$ which is impossible since $\sigma^\A$ and $\sigma^\B$ are optimal policies. The last equality in (ii) follows from (i) since a persistently optimal policy $\sigma^\B$ is optimal.
\end{proof}

\begin{theorem} \label{thmA2}
For every horizon $N=1,2,\ldots$ or $N=\infty$ and for every policy $\pi^\B$ for Nature, there is a nonrandomized policy $\phi^\A$ for the DM such that
\begin{equation}\label{eqthA2}
v_N(x,\phi^\A,\pi^\B)=\min_{\pi^\A\in\Pi^\A} v_N(x,\pi^\A,\pi^\B),\qquad x\in {\bf X}.
\end{equation}
\end{theorem}
\begin{proof}
If Nature plays a policy $\pi^\B,$ the DM has an MDP with the states $h_t=x_0,a_0,b_0,x_1,a_1,b_1,x_2,\ldots,x_t,$ $t=0,1,\ldots.$  Each action set $A(h_t)=A(x_t)$ is finite. For a discounted MDPs with finite action sets there is a Markov optimal policy $\phi^\A$ \cite{Sch75} or \cite[Theorem 3.1]{FK20}. This policy is a nonrandomized policy for the DM, and it satisfies \eqref{eqthA2}.
\end{proof}

The following theorem presents the minimax equation.  We recall that, unlike the case of stochastic games with simultaneous decisions, these equations are not symmetric since Nature knows current decisions of the DM.
\begin{theorem} \label{thmA3}
For every horizon $N=1,2,\ldots$ or $N=\infty,$ 
\begin{equation}\label{eqthmA3}
v_N(x)=\min_{\pi^\A\in\Pi^\A} \sup_{\pi^\B\in\Pi^\B} v_N(x,\pi^\A,\pi^\B)=\max_{\pi^\B\in\Pi^\B} \min_{\pi^\A\in\Pi^\A} v_N(x,\pi^\A,\pi^\B),\qquad x\in {\bf X}.
\end{equation}
\end{theorem}
\begin{proof}
Let us fix $N,$ and let $\pi_*^\A$ and $\pi_*^\B$ be optimal  policies for the DM and Nature whose existence is stated in Theorem~\ref{thmA1}.  Let us fix  $x\in {\bf X}.$ As follows from the definition of optimal policies,
\begin{equation}\label{eqA14}
v_N(x)=\min_{\pi^\A\in\Pi^\A}v_N(x,\pi^\A,\pi_*^\B)\le \inf_{\pi^\A\in\Pi^\A} \sup_{\pi^\B\in\Pi^\B} v_N(x,\pi^\A,\pi^\B)\le  \max_{\pi^\B\in\Pi^\B} v_N(x,\pi_*^\A,\pi^\B)=v_N(x),
\end{equation}
where the inequalities follow from the properties of infimums. Thus, all inequalities in \eqref{eqA14} are equalities, and the infimum in \eqref{eqA14} is the minimum since it is achieved at $\pi^\A=\pi^\A_*.$  The first equality in \eqref{eqthmA3} is proved.  In addition,
\[
v_N(x)=\min_{\pi^\A\in\Pi^\A}v_N(x,\pi^\A,\pi_* ^\B)\le \sup_{\pi^\B\in\Pi^\B} \min_{\pi^\A\in\Pi^\A} v_N(x,\pi^\A,\pi^\B)\le  \inf_{\pi^\A\in\Pi^\A}\sup_{\pi^\B\in\Pi^\B} v_N(x,\pi^\A,\pi^\B)\le v_N(x),\] 
where the  equality follows from the definition of an optimal policy,   the first inequality follows from the properties of supremums and from Theorem~\ref{thmA2}, the second inequality follows from the inequality between $\sup\inf$ and $\inf\sup,$    and the last inequality is taken from \eqref{eqA14}.  Thus, these inequalities hold in the form of equalities, the infimum is the minimum, the first supremum is the maximum, and the second equality in \eqref{eqthmA3} is proved. \end{proof}

The following statement generalizes Theorem~\ref{thmA3}.
\begin{corollary} \label{cor4}
Let $N=1,2,\ldots$ or $N=\infty,$  and let $\Sigma_N^\A\subset\Pi^\A$ and $\Sigma_N^\B\subset\Pi^\B$ be subsets of policies for the DM and Nature such that there  exist an optimal policy $\pi_*^\A\in\Pi^\A$ for the DM and an optimal policy $\pi_*^\B\in\Pi^\B$ for Nature. Then  
\begin{equation*}
v_N(x)=\min_{\pi^\A\in\Sigma_N^\A} \sup_{\pi^\B\in\Sigma_N^\B} v_N(x,\pi^\A,\pi^\B)=\max_{\pi^\B\in\Sigma_N^\B} \min_{\pi^\A\in\Sigma_N^\A} v_N(x,\pi^\A,\pi^\B),\qquad x\in {\bf X}.
\end{equation*}
\end{corollary}
\begin{proof} We observe that, since an optimal policy $\pi_*^\A\in\Sigma_N^\A,$ for $x\in{\bf X}$
\[
v_N(x)=v_N(x,\pi_*^\A,\pi_*^\B)= \min_{\pi^\A\in\Sigma_N^\A} v_N(x,\pi^\A,\pi_*^\B)= \min_{\pi^\A\in\Sigma_N^\A} v_N(x,\pi^\A,\pi_*^\B)=v_N(x).
\]
Therefore, similarly to the proof of Theorem~\ref{thmA3}, for $x\in {\bf X},$
\[
v_N(x)=\min_{\pi^\A\in\Sigma_N^\A}v_N(x,\pi^\A,\pi_*^\B)\le \inf_{\pi^\A\in\Sigma_N^\A} \sup_{\pi^\B\in\Sigma_N^\B} v_N(x,\pi^\A,\pi^\B)\le  \max_{\pi^\B\in\Sigma_N^\B} v_N(x,\pi_*^\A,\pi^\B)=v_N(x),
\]
and
\[
v_N(x)=\min_{\pi^\A\in\Sigma_N^\A}v_N(x,\pi^\A,\pi_* ^\B)\le \sup_{\pi^\B\in\Sigma_N^\B} \min_{\pi^\A\in\Sigma_N^\A} v_N(x,\pi^\A,\pi^\B)\le  \inf_{\pi^\A\in\Sigma_N^\A}\sup_{\pi^\B\in\Sigma_N^\B} v_N(x,\pi^\A,\pi^\B)\le v_N(x).\]
\end{proof}

\end{appendices}
\end{document}